\documentclass[12pt,leqno,twoside]{amsart}
\usepackage{amsmath,amscd,amssymb,amsfonts,latexsym,mathrsfs,bm,extarrows,graphicx,stackrel}
\usepackage[all,cmtip]{xy}
\usepackage[normalem]{ulem}
\usepackage{tikz}

\usepackage{amssymb}
\usepackage{latexsym}
\usepackage{graphicx}
\usepackage{epsfig}
\usepackage{srcltx}
\usepackage[colorlinks=true, pdfstartview=FitV, linkcolor=blue, citecolor=blue, urlcolor=blue]{hyperref}
\usepackage{subcaption}

\newtheorem{theorem}{Theorem}[section]
\newtheorem{corollary}[theorem]{Corollary}
\newtheorem{lemma}[theorem]{Lemma}
\newtheorem{proposition}[theorem]{Proposition}

\theoremstyle{remark}
\newtheorem{remark}[theorem]{\sc Remark}
\theoremstyle{remark}

\theoremstyle{definition}
\newtheorem{definition}[theorem]{Definition}
\theoremstyle{remark}
\newtheorem{example}[theorem]{\sc Example}

\theoremstyle{remark}

\theoremstyle{remark}

\numberwithin{equation}{section}  
\def\be{\begin{equation}}
\def\ee{\end{equation}}
\def\bt{\begin{theorem}}
\def\et{\end{theorem}}
\newcommand{\lra}{\longrightarrow}
\def\bc{\begin{corollary}}
\def\ec{\end{corollary}}
\def\br{\begin{remark}}
\def\er{\end{remark}}
\def\bex{\begin{example}}
\def\eex{\end{example}}
\def\bp{\begin{proposition}}
\def\ep{\end{proposition}}

\usepackage{srcltx}
\usepackage[colorlinks=true, pdfstartview=FitV, linkcolor=blue, citecolor=blue, urlcolor=blue]{hyperref}
\usepackage{subcaption} 
\usepackage{amsbsy}
\usepackage{amsthm,alltt,dsfont}

\usepackage[margin=1.2in]{geometry}

\usepackage[utf8]{inputenc}
\usepackage{amsthm,soul}
\usepackage{enumitem}
\usepackage{setspace,kantlipsum}
\usepackage{tikz-cd}
\usepackage{textcomp}

\usepackage{colonequals}
\usepackage{mathtools,comment}

\usepackage{graphicx,bm}
\usepackage{mathptmx}

\newcommand{\C}{\mathbb{C}}

\newcommand{\Q}{\mathbb{Q}}

\newcommand{\D}{\mathbb{D}}

\newcommand{\cF}{\mathcal{F}}
\newcommand{\cG}{\mathcal{G}}
\newcommand{\cH}{\mathcal{H}}

\newcommand{\cK}{\mathcal{K}}

\newcommand{\cM}{\mathcal{M}}

\newcommand{\cP}{\mathcal{P}}

\newcommand{\cS}{\mathcal{S}}

\newcommand{\cV}{\mathcal{V}}

\newcommand{\Nu}{\mathsf V}

\newcommand{\bH}{\mathbb{H}}

\newcommand{\bP}{\mathbb{P}}

\newcommand{\bZ}{\mathbb{Z}}

\newcommand{\of}{\overline{f}}
\newcommand{\oa}{\overline{a}}
\newcommand{\oX}{\overline{X}}
\newcommand{\oF}{\overline{F}}
\newcommand{\ok}{\overline{k}}
\newcommand{\oi}{\overline{i}}
\newcommand{\oY}{\overline{Y}}

\newcommand{\og}{\overline{g}}
\newcommand{\oG}{\overline{G}}

\newcommand{\oZ}{\overline{Z}}

\newcommand{\Gr}{\mathrm{Gr}}

\newcommand{\MHM}{\mathrm{MHM}}

\newcommand{\Perv}{\operatorname{Perv}}

\newcommand{\im}{\operatorname{im}}
\newcommand{\rk}{\operatorname{rank}}
\newcommand{\Supp}{\operatorname{Supp}}
\newcommand{\id}{\operatorname{id}}

\newcommand{\rat}{{\it {rat}}}
\newcommand{\Hlim}{H_{\mathrm{lim}}}

\begin{document}

\title[Fibers of polynomial maps]
{On the Hodge theory of fibers of complex polynomial maps}

\author{Lauren\c{t}iu Maxim}
\address{L. Maxim : Department of Mathematics, University of Wisconsin-Madison, 480 Lincoln Drive, Madison WI 53706-1388, USA, \newline {\text and} \newline Institute of Mathematics of the Romanian Academy, P.O. Box 1-764, 70700 Bucharest, ROMANIA}
\email{maxim@math.wisc.edu}

\author{John Messina}
\address{J. Messina : Department of Mathematics, University of Wisconsin-Madison, 480 Lincoln Drive, Madison WI 53706-1388, USA}
\email{jamessina@wisc.edu}

\dedicatory{Dedicated to the memory of our friends and colleagues Mihai Tib\u{a}r and Maria Aparecida Soares Ruas.}

\begin{abstract}
We give Hodge-theoretic counterparts and refinements of our
results on the topology of fibers of complex polynomial maps.
We derive upper bounds for the Hodge numbers of vanishing
cohomology and atypical fibers in their first possibly
nonvanishing degrees, in terms of transversal Milnor cohomology
and monodromy along singular strata. We also bound the Hodge
numbers of the non-unipotent part of the limit mixed Hodge
structure on general fiber cohomology and of the image of
nilpotent monodromy on its unipotent part. Finally, we obtain rank bounds for local monodromy at a bifurcation value and monodromy at infinity,
and a Hodge-filtered refinement of lower semicontinuity of the
top Betti number of a general fiber.
\end{abstract}

\date{\today}





\maketitle

\tableofcontents

\section{Introduction}

The study of the Hodge theory of fibers of complex polynomials is an
active area of research in singularity theory, algebraic geometry and topology.
Besides the topology of a fiber, one is interested to understand its
mixed Hodge structure, the action of monodromy, and how these change
when the polynomial varies. Deligne's mixed Hodge theory
\cite{D71,D74}, the theory of limit Hodge structures
\cite{Sch73,Ste76,SZ85}, and Saito's mixed Hodge modules
\cite{Sa88,Sa90} provide the tools for studying these questions.
The mixed Hodge structure on Milnor fiber cohomology and its
compatibility with monodromy are discussed, e.g., in \cite{SS85}, but it can also be captured by mixed Hodge module theory via the nearby/vanishing cycle formalism.
For polynomials, however, singularities in the affine space account
for only part of the vanishing cohomology, and the behavior at infinity
must also be taken into account.

The topology and monodromy of polynomial maps have been studied
extensively, beginning with the case of isolated singularities,
including singularities at infinity, e.g., see
\cite{Br88,ST95,ST00,NN00,DN01,T07}.
Hodge-theoretic methods play an important role in this study. For instance, 
Dimca--Saito \cite{DS93,DS03} relate the cohomology of a general
fiber and its weights to monodromy.
Garc\'ia L\'opez--N\'emethi \cite{GLN96,GLN99a,GLN99b} study
monodromy at infinity and the associated limit mixed Hodge
structure.
Sabbah \cite{Sab97} studies monodromy at infinity by means of the
Fourier transform. Matsui--Takeuchi \cite{MT13} obtain formulas
for Jordan blocks in terms of Newton polyhedra, while
Takeuchi--Tib\u{a}r \cite{TT16} also treat non-tame polynomials
with possibly non-isolated singularities.

In \cite{MM26}, we obtained upper bounds for the first possibly
nonzero Betti numbers of general and atypical fibers in terms of
transversal Milnor numbers. The purpose of this paper is to refine
those bounds by using the mixed Hodge module structures on the
complexes appearing in their proofs. The resulting inequalities
distinguish Hodge types and monodromy eigenvalues, and also take
into account the monodromy of local systems along the singular
strata. They apply without a tameness or Newton nondegeneracy
assumption and allow non-isolated singularities both in the affine
space and at infinity. The Milnor number bounds of \cite{MM26}
are sharp, but for a given polynomial the additional Hodge and
monodromy information may yield smaller upper bounds.

\subsection{Setup}

Let $f:\C^{n+1}\to\C$ be a nonconstant polynomial, with $n\geq1$,
and assume that $\dim\operatorname{Sing}(f)<n$.
Let $B\subset\C$ be its bifurcation set. We denote a general fiber of $f$ 
by $F$, and the atypical fiber over $b\in B$ by $F_b=f^{-1}(b)$.
As in \cite{MM26}, we identify $\C^{n+1}$ with the graph $X$ of
$f$ and consider its closure
$\oX\subset\bP^{n+1}\times\C$, with open affine inclusion
$j:X\hookrightarrow\oX$ and proper projection
$\of:\oX\to\C$. The fiber $\oF_b=\of^{-1}(b)$ is the projective closure
of $F_b$. Let
\[
 \Sigma_b=\Supp\left(\varphi_{\of-b}Rj_*\Q_X\right),
 \quad s_b=\dim\Sigma_b,
 \quad s=\max_{b\in B}s_b,
\]
where we assume $B\ne\emptyset$ when using $s$.
The affine and infinite parts of this support are denoted by
$\Sigma_b^a$ and $\Sigma_b^\infty$, respectively.
We work with algebraic Whitney stratifications adapted to these
complexes, as described in Section~\ref{sec:graphc}.

For a small disc $D_b$ centered at $b\in B$, let $T(F_b)=f^{-1}(D_b)$
be the corresponding tube. The vanishing cohomology at $b$ is defined as
\[
 V^k(b)=H^{k+1}(T(F_b),F;\Q),
\]
with a nearby general fiber $F\subset T(F_b)$.
Its description by vanishing cycles endows it with a mixed Hodge
structure, which we denote by $V_H^k(b)$, see Definition~\ref{def1}.
The vanishing result of \cite[Proposition 3.1]{Br03}, see also \cite[Corollary~2.3]{MM26}, then implies that
\[
 V_H^k(b)=0,\quad\text{for }k\notin[n-s_b,n].
\]
One of our goals in this paper is to study the lowest degree $n-s_b$ when $s_b>0$.

\subsection{Statement of results}
Let
\[
 \cP_b^H:=\varphi_{\of-b}^H
 \left(j_*\Q_X^H[n+1]\right)
\]
be the vanishing cycle mixed Hodge module on $\oF_b$,
with support $\Sigma_b$.
Its restriction to each connected $s_b$-dimensional stratum $S_{b,i}\subset\Sigma_b$ has the form
\[
 \cP_b^H|_{S_{b,i}}\simeq\cV_{b,i}[s_b],
\]
where $\cV_{b,i}$ is an admissible variation of mixed Hodge
structure on $S_{b,i}$.
In Theorem~\ref{thm1} we show that when $s_b>0$ there is a (strict)
monomorphism of mixed Hodge structures
\begin{equation}\label{eq1}
 V_H^{n-s_b}(b)\hookrightarrow
 \bigoplus_i H^0(S_{b,i};\cV_{b,i}),
\end{equation}
compatible with vanishing cycle monodromy.
Recall that the global sections on each stratum identify with the invariants
of its local system monodromy. Thus \eqref{eq1} bounds each
equivariant Hodge number by the corresponding Hodge numbers of
these invariant subspaces.
Moreover, if $s_b\geq2$ and there is no stratum of dimension $s_b-1$ in
$\Sigma_b$, the map \eqref{eq1} is an isomorphism.

The local terms in \eqref{eq1} are described in
Section~\ref{sec5}. For a top-dimensional stratum $S$, let $\Nu_S$ be the stalk
of the variation, and let $MF_S^\pitchfork$ be the transversal
Milnor fiber of $\of-b$ along $S$. Set
\[
 H_S^\pitchfork=\widetilde H^{n-s_b}(MF_S^\pitchfork;\Q),
\]
with its vanishing mixed Hodge structure and Milnor monodromy.
On an affine stratum, one has the identification $\Nu_S\simeq H_S^\pitchfork$ compatibly with monodromy. However, at infinity, duality introduces different Tate twists on the unipotent, resp., non-unipotent summands.
More precisely, with $d_S=n+1-s_b$,
we show in Theorem~\ref{thm2} that there are isomorphisms of mixed Hodge structures
\begin{equation}\label{eq2}
 (\Nu_S)_1\simeq(H_S^\pitchfork)_1^\vee(-d_S),
\end{equation}
\begin{equation}\label{eq2b}
 (\Nu_S)_\lambda\simeq
 (H_S^\pitchfork)_{\lambda^{-1}}^\vee(1-d_S), 
 \quad \lambda\ne1.
\end{equation}
Here, generalized eigenspaces are taken after complexification.
These formulas make use of Saito's duality theorem for vanishing cycles
\cite{Sa89}, see also \cite{MSS20}.
Together with \eqref{eq1}, they yield the explicit affine and
infinity estimates stated in Corollary~\ref{cor2}.

\smallskip

Our next results concern the monodromy $T_b$ of the fibration
$f:T(F_b)\setminus F_b\to D_b\setminus\{b\}$ on the cohomology
of its general fiber $F$.
The nearby cycles of $\of-b$ with coefficients in
$j_*\Q_X^H[n+1]$ define the limit mixed Hodge structure
$\Hlim^k(F;b)$. Its underlying rational vector space identifies
with $H^k(F;\Q)$, with nearby cycle monodromy corresponding
to $T_b$. Its Hodge and weight filtrations describe the degeneration
at $b$ and differ in general from those of Deligne's mixed
Hodge structure on $H^k(F;\Q)$.

In Proposition~\ref{prop3} we compare this limit
mixed Hodge structure with the vanishing mixed Hodge
structure $V_H^k(b)$. For $\lambda\ne1$, we show that the canonical map
identifies their generalized $\lambda$-eigenspaces,
compatibly with the Hodge and weight filtrations.
On the unipotent part of $\Hlim^k(F;b)$, the logarithm of
the monodromy corresponding to $T_b$ defines a morphism
of mixed Hodge structures
\[
 N_b:\Hlim^k(F;b)_1\longrightarrow\Hlim^k(F;b)_1(-1).
\]
Its image, with the mixed Hodge structure induced from the
target, is a subquotient of $V_H^k(b)_1$.
In degree $n-s_b$, Corollary~\ref{cor3} provides upper 
bounds  for the Hodge numbers of the non-unipotent part of the
limit mixed Hodge structure and of $\operatorname{im}N_b$.

On underlying rational vector spaces, we also have
\[
 \operatorname{rk}(T_b-\operatorname{id})\leq\dim_\Q V^k(b).
\]
Combining this rank inequality with the bounds of
\cite{MM26}, we obtain bounds for
$\operatorname{rk}(T_b-\operatorname{id})$ and for the number
and sizes of nontrivial Jordan blocks with eigenvalue $1$
in degree $n-s$.
By expressing monodromy at infinity as a product of local
monodromies, we then obtain a rank bound there as well
(Corollary~\ref{cor6}).
For more results on the eigenvalues and Jordan blocks
of monodromy at infinity see, e.g., 
\cite{GLN96,GLN99a,Sab97,DS03,MT13,TT16}.

\smallskip

For an atypical fiber, in \cite[Theorem~4.8]{MM26} we gave upper bounds for the first
possibly nonzero reduced cohomology group, in degree $n-s-1$,
using only strata at infinity.
In Theorem~\ref{thm3} of this paper we refine that result as follows.
For a connected $s$-dimensional stratum
$S\subset\Sigma_b^\infty$, set
\[
 K_S=H^0(S;\cV_{S,1}),\quad
 N_S:K_S\longrightarrow K_S(-1),
\]
where the subscript $1$ refers to the unipotent part of vanishing cycle monodromy
and $N_S$ is induced by its nilpotent logarithm.
Then we show that there is a monomorphism
\begin{equation}\label{eq3}
 \widetilde H^{n-s-1}(F_b;\Q)
 \hookrightarrow
 \bigoplus_{\substack{S\subset\Sigma_b^\infty\\\dim S=s}}
 \ker\left(N_S(1):K_S(1)\longrightarrow K_S\right).
\end{equation}
In particular, as already observed in \cite{Br03,MM26}, if $\dim\Sigma_b^\infty<s$ the group on the left vanishes.
The twist in \eqref{eq3} comes from the variation morphism.
Combining this injection with \eqref{eq2} implies
\begin{equation}\label{eq4}
 h^{p,q}\left(\widetilde H^{n-s-1}(F_b;\Q)\right)
 \leq
 \sum_{\substack{S\subset\Sigma_b^\infty\\\dim S=s}}
 h_1^{n-s-p,n-s-q}(H_S^\pitchfork).
\end{equation}
Let us note that the intermediate bound in \eqref{eq3} can be smaller:
it uses invariance under monodromy along $S$ and the kernel
of $N_S$, in addition to the eigenvalue-$1$ condition.
The extension of Broughton's example in
Example~\ref{ex1} shows that the Hodge bound for the 
atypical fiber can be an equality, even when
the monodromy rank bound is strict.

\smallskip

Finally, we consider a one-parameter family $f_\tau=P(-,\tau)$
of degree $d$ polynomials and a common general value $c$, with
$F^\tau=f_\tau^{-1}(c)$.
The sheaf-theoretic proof of lower semicontinuity of the top
Betti number in \cite[Theorem~5.2]{MM26} admits a mixed Hodge
module version, which we discuss in Theorem~\ref{thm4}.
More precisely, if $H_{c,\mathrm{lim}}^n(F^\tau;0)$ denotes the limit mixed
Hodge structure on the compactly supported cohomology
$H_c^n(F^\tau;\Q)$ as $\tau$ approaches $0$, then we show in 
Theorem~\ref{thm4}  that there is a monomorphism
of mixed Hodge structures
\[
 H_c^n(F^0;\Q)\hookrightarrow
 H_{c,\mathrm{lim}}^n(F^\tau;0).
\]
As the Hodge-graded dimensions of the limit agree with those
of a nearby fiber,  strictness and Poincar\'e duality then imply,
for each $p$ and sufficiently small $\tau\ne0$, that
\begin{equation}\label{eq5}
 \dim_\C\Gr_F^p H^n(F^0;\C)
 \leq\dim_\C\Gr_F^p H^n(F^\tau;\C),
\end{equation}
with the same inequality also holding for compactly supported cohomology.
Summing over $p$ recovers the lower semicontinuity of the
top Betti number proved in \cite[Theorem~5.2]{MM26}, see also \cite[Proposition 2.1]{ST07} for a proof by homotopy-theoretic methods.
A related semicontinuity result for the spectrum at infinity
of analytic families of weakly tame regular functions was
proved by N\'emethi--Sabbah \cite{NS99}. As an application of our semicontinuity result \eqref{eq5}, in Proposition \ref{prop4} we obtain sharp upper bounds for the dimensions of Hodge graded pieces of top cohomology of the generic fiber of a polynomial $f:\C^{n+1} \to \C$.

\smallskip

The proofs of our results follow the stratified arguments of \cite{MM26},
using mixed Hodge modules to obtain morphisms compatible
with the Hodge and weight filtrations. Strictness turns these
morphisms into Hodge-number inequalities, while duality and
the variation morphism determine the Tate twists in the
contributions from infinity.

\subsection{Acknowledgments} 
We would like to dedicate this work to our colleagues, friends and mentors Mihai Tib\u{a}r and Maria Aparecida Soares Ruas, 
whose significant contributions to mathematics, generosity, wisdom, and lasting influence on our community will continue to inspire us and future generations.
L. Maxim acknowledges support from the Simons Foundation (award MPS-TSM-00007095). 


\section{Mixed Hodge structures, monodromy, nearby and vanishing cycles}

This section fixes notation and conventions, and recalls relevant background for the standard theory used in this paper. References for the material presented here include \cite{D71,D74,Sa88,Sa90,D04,M19,MS22}.

\subsection{Mixed Hodge structures and Tate twists}\label{sec21}

A mixed Hodge structure (MHS, for short) is a finite dimensional $\Q$-vector
space $H$ endowed with a finite increasing weight filtration $W_\centerdot$, and with a finite decreasing
Hodge filtration $F^\centerdot$ on $H_\C=H\otimes \C$,  such that $\left(\Gr_k^W H, F^\centerdot\right)$ is a pure Hodge structure of weight $k$, for all $k$.
A morphism of mixed Hodge structures is a $\Q$-linear map between two mixed Hodge structures, which is compatible with both filtrations $F^\centerdot$ and  $W_\centerdot$.

The Hodge numbers of a mixed Hodge structure $H$ are defined as
\[
 h^{p,q}(H)=
 \dim_\C\Gr_F^p\Gr_{p+q}^W(H_\C).
\]
Our convention is that $\Q(1)$ has type $(-1,-1)$ and weight $-2$. Hence,
\begin{equation}\label{eq6}
 h^{p,q}\left(H(-r)\right)=h^{p-r,q-r}(H) \ \ {\rm and} \ \ 
 h^{p,q}\left(H^\vee(-r)\right)=h^{r-p,r-q}(H).
\end{equation}

A morphism of mixed Hodge structures is strict for both $F$ and
$W$.  Therefore, if $K\hookrightarrow H$ is a monomorphism of mixed Hodge structures, then
\begin{equation}\label{eq7}
 h^{p,q}(K)\le h^{p,q}(H)
\end{equation}
for all $p,q$.

Assume next that $H$ carries a monodromic mixed Hodge structure with
quasi-unipotent monodromy $T=T_sT_u$.  Thus $T_s$ is a finite-order
automorphism of the mixed Hodge structure, while
\[
 N=\log T_u:H\longrightarrow H(-1)
\]
is a morphism of mixed Hodge structures.  An eigenvalue $\lambda$ of $T_s$ is a root of unity, and we let
\[
 H_\lambda=\ker(T_s-\lambda)\subset H_\C,
\]
that is, the generalized $\lambda$-eigenspace of $T$.  
We set
\[
 h_\lambda^{p,q}(H)=
 \dim_\C\left(\Gr_F^p\Gr_{p+q}^W H_\C\right)_\lambda,
\]
for the dimension of the $\lambda$-eigenspace of $T_s$ on $\Gr_F^p\Gr_{p+q}^W H_\C$. 
For a monomorphism $K\hookrightarrow H$ compatible with monodromy we then have
\begin{equation}\label{eq8}
 h_\lambda^{p,q}(K)\le h_\lambda^{p,q}(H).
\end{equation}

\subsection{Mixed Hodge modules}
For a complex algebraic variety $Z$, let $\MHM(Z)$ denote Saito's
category of graded-polarizable mixed Hodge modules \cite{Sa90}, and let
\[
 \rat:D^b\MHM(Z)\longrightarrow D_c^b(Z;\Q)
\]
be the rational forgetful functor, assigning the underlying rational constructible complexes. The functor $\rat$ is faithful and exact
on $\MHM(Z)$, with values in the category $\Perv(Z;\Q)$ of perverse sheaves on $Z$.  The standard operations
$f_*,f_!,f^*,f^!$, Verdier duality $\D$, and nearby and vanishing cycles
exist in $D^b\MHM$ and are compatible with their constructible sheaf
counterparts.

If $a:Z\to\mathrm{pt}$ is the constant map to a point space and
$M\in D^b\MHM(Z)$, the objects
$H^k(a_*M)$ and $H^k(a_!M)$
are mixed Hodge structures.  All maps between such groups obtained from the functorial
operations are morphisms of mixed Hodge structures and are therefore
strict.

If $Z$ is smooth of pure dimension $d$, then
$\Q_Z^H[d]\in\MHM(Z)$ is pure of weight $d$, with underlying rational complex 
 $\Q_Z[d]\in \Perv(Z;\Q)$, and
\begin{equation}\label{eq9}
 \D\Q_Z^H[d]\simeq\Q_Z^H[d](d).
\end{equation}
If $j:Z\hookrightarrow\oZ$ is an affine open embedding, both $j_!$
and $j_*$ are $t$-exact for the standard $t$-structures on the derived
categories of mixed Hodge modules.

\subsection{Nearby and vanishing cycles}
Let $g:Z\to\C$ be a regular function and let $M\in\MHM(Z)$.
The shifted nearby and vanishing cycle functors $\psi_g[-1]$,
resp., $\varphi_g[-1]$, preserve perverse sheaves and they lift
to functors $\psi_g^H$ and $\varphi_g^H$ that preserve mixed
Hodge modules.
Their underlying rational complexes therefore satisfy
\[
 \rat(\psi_g^H M)=\psi_g[-1](\rat(M)),
 \quad
 \rat(\varphi_g^H M)=\varphi_g[-1](\rat(M)).
\]

The underlying rational complexes of $\psi_g^H M$ and $\varphi_g^H M$
carry canonical quasi-unipotent monodromy automorphisms,
induced by a positively oriented loop around $0$ in the
target $\C$ of $g$. The decompositions of these rational complexes
into unipotent and non-unipotent parts lift to direct-sum
decompositions in $\MHM(g^{-1}(0))$, namely
\[
 \psi_g^H M
 =\psi_{g,1}^H M\oplus\psi_{g,\ne1}^H M,
 \quad
 \varphi_g^H M
 =\varphi_{g,1}^H M\oplus\varphi_{g,\ne1}^H M.
\]
The subscript $1$ denotes the generalized eigenvalue-$1$
summand, on which monodromy is unipotent, and $\ne1$ denotes
the sum of the remaining generalized eigenspaces.
Individual generalized $\lambda$-eigenspaces are considered
after complexification, and they are not, in general, rational
mixed Hodge modules.

On the unipotent summands there are morphisms
\begin{equation}\label{eq10}
 \operatorname{can}:\psi_{g,1}^H M\rightarrow\varphi_{g,1}^H M,
 \quad
 \operatorname{Var}:\varphi_{g,1}^H M
 \rightarrow\psi_{g,1}^H M(-1).
\end{equation}
These canonical and variation morphisms satisfy
\begin{equation}\label{eq11}
\begin{aligned}
 \operatorname{Var}\circ\operatorname{can}
 &=N:\psi_{g,1}^H M\longrightarrow\psi_{g,1}^H M(-1),\\
 \operatorname{can}(-1)\circ\operatorname{Var}
 &=N:\varphi_{g,1}^H M\longrightarrow\varphi_{g,1}^H M(-1).
\end{aligned}
\end{equation}
Here $N$ denotes the logarithm of unipotent monodromy.
On the non-unipotent parts, the canonical morphism is an
isomorphism in $\MHM(g^{-1}(0))$:
\[
 \operatorname{can}:\psi_{g,\ne1}^H M
 \overset{\sim}{\lra} \varphi_{g,\ne1}^H M.
\]
After complexification, $\operatorname{can}$ identifies the corresponding
generalized $\lambda$-summands for every $\lambda\ne1$.

The following duality isomorphisms
follow from Saito's duality theorem for vanishing cycles
\cite{Sa89}, see also
\cite[(4.5.1) and \S4.7]{MSS20}:
\begin{align}
 \varphi_{g,1}^H\D M
 &\simeq\D\varphi_{g,1}^H M,\label{eq12}\\
 \psi_{g,\lambda}^H(\D M)
 &\simeq
 \left(\D(\psi_{g,\lambda^{-1}}^H M)\right)(1),
 \quad\lambda\ne1.\label{eq13}
\end{align}
The inverse eigenvalue in \eqref{eq13} comes from dual
monodromy.

\subsection{Global sections of admissible variations}

Let $S$ be a smooth and connected algebraic variety, and let $\cV$ be a
graded-polarizable admissible variation of mixed Hodge structure on
$S$. The group $H^0(S;\cV)$
is a mixed Hodge structure, and evaluation at $x\in S$ is a strict
monomorphism
\begin{equation}\label{eq14}
 H^0(S;\cV)\hookrightarrow\cV_x.
\end{equation}
On underlying vector spaces, its image is $\cV_x^{\pi_1(S,x)}$.
Thus $H^0(S;\cV)$ can be identified with the subspace invariant under
local system monodromy, and it comes equipped with the Hodge and weight
filtrations induced from $\cV_x$.

If the underlying local system also carries a quasi-unipotent
automorphism, evaluation commutes with this automorphism.
After complexification, the evaluation map respects the generalized
eigenspace decomposition and commutes with the logarithm of the
unipotent part. For the vanishing cycle variations below,
this additional automorphism is the canonical vanishing cycle monodromy, while the local system monodromy will be referred to as the vertical monodromy along the stratum.


\section{Vanishing mixed Hodge structure at a bifurcation value}

\subsection{The graph completion}\label{sec:graphc}
Let $f:\C^{n+1}\to\C$ be a nonconstant polynomial, with $n\geq1$,
and assume throughout that
\[
 \dim\operatorname{Sing}(f)<n.
\]
Let $B\subset\C$ be the bifurcation set of $f$, that is, the finite set of
values away from which $f$ is a locally trivial topological fibration.
For $b\in B$, let
\[
 F_b=f^{-1}(b)
\]
be the corresponding atypical fiber.
The fibers over $\C\setminus B$ are the general fibers of $f$. 
They are smooth of dimension $n$ and are connected under the above 
assumption on $\dim\operatorname{Sing}(f)$. We will denote a general fiber by $F$.

Set $d=\deg f$ and identify $\C^{n+1}$ with the graph $X$ of $f$
by the map $x\mapsto(x,f(x))$. Its closure in
$\bP^{n+1}\times\C$ is
\begin{equation}\label{eq15}
 \oX=\left\{([x_0:x],t)\ \middle|\
 f^h(x_0,x)-t x_0^d=0\right\},
\end{equation}
where $f^h$ is the degree-$d$ homogenization of $f$ in the variable $x_0$.
Let
\[
 j:X\hookrightarrow\oX,\quad
 \of:\oX\longrightarrow\C
\]
be the open inclusion and the second projection.
The map $\of$ is proper and restricts to $f$ on $X$.
Its fiber over $b$ is the projective closure of $F_b$.
We denote
\[
 X^\infty=\oX\setminus X,\quad
 \oF_b=\of^{-1}(b),\quad
 F_b^\infty=\oF_b\cap X^\infty.
\]

Since $X$ is smooth of dimension $n+1$ and $j$ is an affine
open inclusion, we have
\begin{equation}\label{eq16}
 \cM:=j_*\Q_X^H[n+1]\in\MHM(\oX).
\end{equation}
For $b\in B$, consider its vanishing cycles along the fiber
$\oF_b$, namely
\begin{equation}\label{eq17}
 \cP_b^H:=
 \varphi_{\of-b}^H\cM
 \in\MHM(\oF_b),
 \quad \ {\rm with} \ \ 
 \Sigma_b=\Supp\rat(\cP_b^H),
\end{equation}
and note that
\begin{equation}\label{eq18}
 \rat(\cP_b^H)
 =\varphi_{\of-b}[-1]\left(Rj_*\Q_X[n+1]\right)
 =\varphi_{\of-b}\left(Rj_*\Q_X\right)[n]
 =:\cP_b.
\end{equation}
We also regard $\cP_b^H$ as a mixed Hodge module on its support
$\Sigma_b$, and use the same convention for $\cP_b$.
The support $\Sigma_b$ is projective, since it is closed in
the projective fiber $\oF_b$.

Let us fix an algebraic Whitney stratification $\cS$ of $\oX$
for which $Rj_*\Q_X$ is constructible, and $X^\infty$ and
$\oF_b$, for all $b\in B$, are unions of strata.
For each $b\in B$, choose an algebraic Whitney stratification
$\cS_b$ of $\oF_b$ refining the restriction of $\cS$,
such that $\Sigma_b$ is a union of strata and $\cP_b$ is
$\cS_b$-constructible.

We define the affine and infinite parts of $\Sigma_b$ as
\[
 \Sigma_b^a=\Sigma_b\cap X
 =F_b\cap\operatorname{Sing}(f) \ \ \ {\rm and} \ \ \ 
 \Sigma_b^\infty=\Sigma_b\cap X^\infty.
\]
On $X$, the complex $Rj_*\Q_X$ restricts to $\Q_X$, so its
vanishing cycles are the usual vanishing cycles of $f$.
Their support on $F_b$ is precisely
$F_b\cap\operatorname{Sing}(f)=\Sigma_b^a$. Note that, since $\dim\operatorname{Sing}(f)<n$, every fiber $F_b$ is a reduced hypersurface, so we can further identify 
\[
 \Sigma_b^a=\operatorname{Sing}(F_b).
\]

\subsection{The vanishing mixed Hodge structure}
In this section we introduce Hodge-theoretic analogues of the
vanishing cohomology at a bifurcation value, which played an
essential role in obtained the bounds for the Betti numbers of polynomial
fibers in \cite{MM26}. 

Let $b\in B$ and set $s_b=\dim\Sigma_b$. If $D_b \subset \C$ is a sufficiently small disc centered at $b$, we denote by $T(F_b):=f^{-1}(D_b)$ the corresponding tube around the atypical fiber $F_b$. Recall from \cite{MM26} that the \emph{vanishing cohomology of $f$ at $b$} is defined as
\[V^k(b):=H^{k+1}(T(F_b),F;\Q),\]
for $F$ a general fiber contained in $T(F_b)$. 

\begin{definition}\label{def1}
Let $a_b:\Sigma_b\to\mathrm{pt}$ be the constant map. The
\emph{vanishing mixed Hodge structure in degree $k$ at $b$} is defined as
\begin{equation}\label{eq19}
 V_H^k(b):=
 H^{k-n}\left((a_b)_*\cP_b^H\right).
\end{equation}
\end{definition}

\begin{proposition}\label{prop1}
The rational vector space underlying $V_H^k(b)$ is 
isomorphic to $V^k(b)$. 
Moreover,
\begin{equation}\label{eq20}
 V_H^k(b)=0\quad {\rm for } \ \ k\notin[n-s_b,n].
\end{equation}
\end{proposition}

\begin{proof}
By \eqref{eq18}, the underlying rational complex of $\cP_b^H$ is
$\varphi_{\of-b}Rj_*\Q_X[n]$, supported on $\Sigma_b$.
Together with \cite[Proposition~2.2]{MM26}, we obtain:
\begin{align*}
 \rat\left(V_H^k(b)\right)
 &\simeq\bH^{k-n}\!\left(
   \Sigma_b;\rat(\cP_b^H)\right)\\
 &\simeq\bH^{k-n}\!\left(
   \oF_b;\varphi_{\of-b}Rj_*\Q_X[n]\right)\\
 &\simeq\bH^k\!\left(
   \oF_b;\varphi_{\of-b}Rj_*\Q_X\right)\\
 &\simeq H^{k+1}(T(F_b),F;\Q)=V^k(b).
\end{align*}

By \cite[Proposition 3.1]{Br03} or \cite[Corollary 2.3]{MM26},
$V^k(b)=0$ for $k\notin[n-s_b,n]$.
So, since $\rat$ is faithful, the same vanishing holds for
$V_H^k(b)$.
\end{proof}


\section{Top-dimensional singular strata and the lowest vanishing group}

Let $\cS_b$ be the Whitney stratification of $\oF_b$ introduced in 
Section~\ref{sec:graphc}, and let $S_{b,i}$ be the connected
$s_b$-dimensional strata contained in $\Sigma_b$.
Since $\rat(\cP_b^H)$ is $\cS_b$-constructible and perverse,
and $S_{b,i}$ has maximal dimension in its support, the restriction
of $\cP_b^H$ to $S_{b,i}$ is an admissible variation of mixed Hodge
structure shifted by $s_b$. Denote this variation by $\cV_{b,i}$
and its underlying local system by $L_{b,i}$. Then
\begin{equation}\label{eq21}
 \cP_b^H|_{S_{b,i}}\simeq\cV_{b,i}[s_b].
\end{equation}
The rational complex $\rat(\cP_b^H)$ carries the 
vanishing cycle monodromy $M_b$, and all restriction and attaching
morphisms below commute with this action.  Moreover, its restriction to
$L_{b,i}$ is quasi-unipotent. For $x \in S_{b,i}$, we denote the induced automorphism by
\[
 M_{b,i}:(L_{b,i})_x
 \longrightarrow(L_{b,i})_x.
\]
Besides the vanishing cycle monodromy $M_{b,i}$, the stalk
$(L_{b,i})_x$ carries the monodromy action of the local system
$L_{b,i}$ along $S_{b,i}$, which we denote by
\[
 A_b^i:\pi_1(S_{b,i},x)\longrightarrow
 \operatorname{Aut}\left((L_{b,i})_x\right).
\]
If $\pi_1(S_{b,i},x)\simeq\bZ$ and a generator is chosen, the resulting
automorphism is the vertical monodromy along $S_{b,i}$ in the
terminology of \cite[\S3.1]{MPT22}.  The image of $A_b^i$
commutes with $M_{b,i}$.  Thus $A_b^i$ records monodromy along the
stratum, whereas $M_{b,i}$ accounts for the  vanishing cycle monodromy.  For an affine stratum, $M_{b,i}$ agrees with the Milnor monodromy on a normal slice.  Its relation to the corresponding transversal Milnor monodromy at infinity is given in
Theorem~\ref{thm2} below.

\medskip

In order to formulate the main result of this section we need to introduce some notation. Following \cite[\S 2]{MM26}, let
\[
 \Sigma_b^\ell=
 \bigcup_{\substack{S\in\cS_b,\ S\subset\Sigma_b\\
                     \dim S=\ell}}S
                 \ \    \ \ \ {\rm and} \ \ \ \ \ 
 U_b^\ell=
 \bigcup_{\substack{S\in\cS_b,\ S\subset\Sigma_b\\
                     \dim S\ge\ell}}S.
\]
Then each $U_b^\ell$ is open in $\Sigma_b$, and
\[
  \Sigma_b^{s_b}=U_b^{s_b}\subset U_b^{s_b-1}\subset\cdots\subset U_b^0=\Sigma_b.
\]
Set
\[
 \cP_b^{H,\ell}:=\cP_b^H|_{U_b^\ell}
 \in\MHM(U_b^\ell),\]
with 
$ \cP_b^{H,0}=\cP_b^H$. 
For the inclusions
\[
 v_\ell:\Sigma_b^\ell\hookrightarrow U_b^\ell,\quad
 u_\ell:U_b^{\ell+1}\hookrightarrow U_b^\ell,
\]
we have
\[
 (u_\ell)^*\cP_b^{H,\ell}\simeq\cP_b^{H,\ell+1},
\]
and there is an attaching triangle in $D^b\MHM(U_b^\ell)$:
\begin{equation}\label{eq22}
 (v_\ell)_!(v_\ell)^!\cP_b^{H,\ell}
 \longrightarrow\cP_b^{H,\ell}
 \longrightarrow
 (u_\ell)_*(u_\ell)^*\cP_b^{H,\ell}
 \overset{[1]}{\lra}.
\end{equation}

\begin{lemma}\label{lem1}
For every $0\le\ell\le s_b$ and every $r<-\ell$, we have
\[
 H^r\!\left(
 (a_{\Sigma_b^\ell})_*
 (v_\ell)^!\cP_b^{H,\ell}\right)=0.
\]
In particular, if $0\le j\le s_b-1$ and
$0\le\ell\le s_b-j-1$, this vanishing holds in degrees
$-s_b+j$ and $-s_b+j-1$.
\end{lemma}

\begin{proof}
Let $K_\ell=\rat(\cP_b^{H,\ell})$, which is a perverse sheaf on
$U_b^\ell$.  Every stratum contained in $\Sigma_b^\ell$ has complex
dimension $\ell$.  Because $K_\ell$ is perverse,  it follows by the costalk condition that
\[
 \cH^q\left((v_\ell)^!K_\ell\right)=0, \ \ \ q<-\ell.
\]

Now consider the hypercohomology spectral sequence
\[
 E_2^{p,q}=
 H^p\!\left(\Sigma_b^\ell;
 \cH^q\left((v_\ell)^!K_\ell\right)\right)
 \Longrightarrow
 \bH^{p+q}\!\left(\Sigma_b^\ell;(v_\ell)^!K_\ell\right).
\]
By the above vanishing, possible nonzero $E_2$ terms have $q\ge-\ell$ and, by sheaf cohomology, $p\ge0$.  Thus every possibly nonzero $E_2$ term has
$p+q\ge-\ell$.  Therefore the hypercohomology in the abutment is zero in every degree
$r<-\ell$.  This proves the claimed vanishing after applying $\rat$, and hence also  in $\MHM(\mathrm{pt})$, since $\rat$ is faithful and exact.  Finally, if
$\ell\le s_b-j-1$, then both $-s_b+j$ and $-s_b+j-1$ are strictly
smaller than $-\ell$.
\end{proof}


The main result of this section is the following Hodge-theoretic generalization of \cite[Theorem 2.12]{MM26}.
\begin{theorem}\label{thm1}
Assume $s_b>0$.  There is a monomorphism of mixed Hodge
structures, compatible with vanishing cycle monodromy,
\begin{equation}\label{eq23}
 V_H^{n-s_b}(b)\hookrightarrow
 \bigoplus_i H^0(S_{b,i};\cV_{b,i}).
\end{equation}
Evaluation at $x \in S_{b,i}$  identifies the underlying vector space of the
$i$-th summand with
\[
 \left((\cV_{b,i})_x\right)^{A_b^i}.
\]
If $s_b\ge2$ and $\Sigma_b$ has no stratum of dimension $s_b-1$,
then \eqref{eq23} is an isomorphism.
\end{theorem}

\begin{proof}
For $0\leq\ell\leq s_b-1$, the long exact sequence of
\eqref{eq22} contains the terms
\[
 \lra H^{-s_b}\!\left(
  \Sigma_b^\ell;(v_\ell)^!\cP_b^{H,\ell}\right)
 \longrightarrow
 H^{-s_b}\!\left(U_b^\ell;\cP_b^{H,\ell}\right)
 \overset{\rho_\ell}{\lra}
 H^{-s_b}\!\left(U_b^{\ell+1};\cP_b^{H,\ell+1}\right) \lra.
\]
Because $-s_b<-\ell$, we get by Lemma~\ref{lem1} that
\[
 H^{-s_b}\!\left(
  \Sigma_b^\ell;(v_\ell)^!\cP_b^{H,\ell}\right)=0,
\]
hence  $\rho_\ell$ is a monomorphism.  Composing these monomorphisms, we get
\begin{align*}
 V_H^{n-s_b}(b)
 =H^{-s_b}\!\left(\Sigma_b;\cP_b^H\right)
 \hookrightarrow
 H^{-s_b}\!\left(U_b^{s_b};\cP_b^{H,s_b}\right)
 \overset{\eqref{eq21}}{\simeq}
 \bigoplus_i H^0(S_{b,i};\cV_{b,i}),
\end{align*}
proving \eqref{eq23}.  Moreover, for $x\in S_{b,i}$, we have by
\eqref{eq14} that 
\[ H^0(S_{b,i};\cV_{b,i}) \hookrightarrow(\cV_{b,i})_x,\]
 with 
\[ \rat \left( H^0(S_{b,i};\cV_{b,i}) \right) =((L_{b,i})_x)^{A_b^i}.\]

Assume now $s_b\geq2$.  If $0\leq\ell\leq s_b-2$, then also
$-s_b+1<-\ell$, so Lemma~\ref{lem1} implies
\[
 H^{-s_b+1}\!\left(
  \Sigma_b^\ell;(v_\ell)^!\cP_b^{H,\ell}\right)=0.
\]
The long exact sequence of \eqref{eq22} therefore makes $\rho_\ell$ an isomorphism for $0\leq\ell\leq s_b-2$.  If
$\Sigma_b^{s_b-1}=\emptyset$, then
$U_b^{s_b-1}=U_b^{s_b}$, and the above composition is an isomorphism.
Finally, note that all arrows are morphisms in $\MHM(\mathrm{pt})$, so in particular they are strict.
\end{proof}

\begin{corollary}\label{cor1}
For every $p,q$ and every eigenvalue $\lambda$ of vanishing cycle
monodromy, we have
\begin{equation}\label{eq24}
 h_\lambda^{p,q}\left(V_H^{n-s_b}(b)\right)
 \le
 \sum_i
 h_\lambda^{p,q}\left(H^0(S_{b,i};\cV_{b,i})\right)
 \le
 \sum_i h_\lambda^{p,q}\left((\cV_{b,i})_x\right).
\end{equation}
\end{corollary}
\begin{proof}
For every $\lambda$, Theorem~\ref{thm1} and
\eqref{eq14} yield strict monomorphisms
\[
 V_H^{n-s_b}(b)_\lambda
 \hookrightarrow
 \bigoplus_i H^0(S_{b,i};\cV_{b,i})_\lambda
 \hookrightarrow
 \bigoplus_i(\cV_{b,i})_{x,\lambda}.
\]
Applying $\Gr_F^p\Gr_{p+q}^W$ implies
\[
 h_\lambda^{p,q}\left(V_H^{n-s_b}(b)\right)
 \leq\sum_i h_\lambda^{p,q}
 \left(H^0(S_{b,i};\cV_{b,i})\right)
 \leq\sum_i h_\lambda^{p,q}\left((\cV_{b,i})_x\right).
\]
\end{proof}


\section{Transversal mixed Hodge structures}\label{sec5}
Let $S=S_{b,i}\subset\Sigma_b$ be a connected top-dimensional stratum,
with
\[
 \cV_S:=\cV_{b,i},\quad L_S:=L_{b,i},\quad M_S:=M_{b,i}.
\]
Thus $\cV_S$ is the admissible variation defined by
\[
 \cP_b^H|_S\simeq\cV_S[s_b],
\]
$L_S$ is its underlying local system, and $M_S$ is the
vanishing cycle monodromy. For $x\in S$, we denote the
stalk at $x$ of the variation $\cV_S$ by
\[
 \Nu_S:=(\cV_S)_x.
\]
Let $i_x:\{x\}\hookrightarrow\oF_b$ be the inclusion, and set
\[
d_S=n+1-s_b.
\]

As in \cite[Lemma~2.10]{MM26}, we choose a general projective
linear subspace $\overline E\subset\bP^{n+1}$ of codimension
$s_b$ meeting $S$ transversally at $x$. Let
$E=\overline E\cap\C^{n+1}$ and $g=f|_E$. Denote the graph of $g$
by $Y$, its closure in $\oX$ by $\oY$, and the restriction
$\of|_{\oY}$ by $\og$. The corresponding inclusions form the Cartesian
square
\[
 \begin{array}{ccc}
 Y&\xrightarrow{\ i\ }&X\\
 {\scriptstyle j_E}\downarrow&&\downarrow{\scriptstyle j}\\
 \oY&\xrightarrow{\ \oi\ }&\oX.
 \end{array}
\]
On the fiber over $b$, let
\[
 \oi_b:\oG_b:=\og^{-1}(b)
 =\oY\cap\oF_b\hookrightarrow\oF_b
\]
be the inclusion. For this general transversal slice, we get as in \cite[Lemma~2.10]{MM26} that
\begin{align}
 \oi_b^*\cP_b^H[-s_b]
 &=\oi_b^*\varphi_{\of-b}^H
   \left(j_*\Q_X^H[n+1]\right)[-s_b]\notag\\
 &\overset{\mathrm{(a)}}{\simeq}
 \varphi_{\og-b}^H
 \left(\oi^*j_*\Q_X^H[n+1]\right)[-s_b]\notag\\
 &\overset{\mathrm{(b)}}{\simeq}
 \varphi_{\og-b}^H
 \left(j_{E*}i^*\Q_X^H[n+1]\right)[-s_b]\notag\\
 &\simeq
 \varphi_{\og-b}^H\left(j_{E*}\Q_Y^H[d_S]\right),
 \label{eq25}
\end{align}
where $\mathrm{(a)}$ follows by base change for vanishing cycles
\cite[Lemma~4.3.4]{Schu03}, and $\mathrm{(b)}$ by base
change for the Cartesian square
\cite[Proposition~4.3.1]{Schu03}.

The support of
$\varphi_{\og-b}^H(j_{E*}\Q_Y^H[d_S])$ is isolated at $x$. 
And, as in \cite[Lemma~2.10]{MM26}, its stalk cohomology at $x$ is related to the cohomology of the transversal Milnor fiber
$MF_S^\pitchfork$ of $\og-b$ at $x$, whose reduced cohomology
is concentrated in degree $n-s_b$. This correspondence will be explained in detail below.

Let
\begin{equation}\label{eq26}
 H_S^\pitchfork:=\widetilde H^{n-s_b}(MF_S^\pitchfork;\Q),
 \quad \mu_S^\pitchfork=\dim H_S^\pitchfork,
\end{equation}
with its vanishing mixed Hodge structure and transversal Milnor
monodromy $h_S$.
With the eigenspace convention of Section~\ref{sec21}, after extending
scalars to $\C$, and for $r\gg0$, we define
\[
 (\Nu_S)_\lambda:=\ker(M_S-\lambda)^r,
 \quad
 (H_S^\pitchfork)_\lambda:=\ker(h_S-\lambda)^r.
\]
Since $M_S$ commutes with the local system monodromy along $S$,
the spaces $(\Nu_S)_\lambda$ are the stalks of generalized eigenspace
subvariations, which we denote by $\cV_{S,\lambda}$.
Under the affine identification of Proposition~\ref{prop2} below,
$M_S$ corresponds to $h_S$. At infinity,
Theorem~\ref{thm2} identifies $M_S$ with the inverse
dual of $h_S$, with corresponding eigenvalues.

\subsection{Affine strata}
\begin{proposition}\label{prop2}
For a top-dimensional affine stratum $S\subset\Sigma_b^a$, there is an
isomorphism of monodromic mixed Hodge structures
\[
 \Nu_S\simeq H_S^\pitchfork
\]
which identifies $M_S$ with $h_S$.  Equivalently, after extending scalars to $\C$, we have for every eigenvalue $\lambda$,
\begin{equation}\label{eq27}
 (\Nu_S)_\lambda\simeq(H_S^\pitchfork)_\lambda.
\end{equation}
\end{proposition}

\begin{proof}
We apply \cite[equation~(2.5)]{MM26} to the transverse
restriction $g=f|_E$. Since $x\in Y$, the base-change
isomorphism \eqref{eq25} yields the identification
\[
 \Nu_S
 \simeq H^0i_x^*\varphi_{g-b}^H\left(\Q_Y^H[d_S]\right)
 \simeq \widetilde H^{n-s_b}(MF_{g-b,x};\Q)=\widetilde H^{n-s_b}(MF_S^\pitchfork;\Q)
 =H_S^\pitchfork.
\]
These are isomorphisms of mixed Hodge structures. Indeed, the first
comes from base change in mixed Hodge modules, and the second
is the vanishing cycle description of the mixed Hodge structure on Milnor fiber cohomology. Both commute with monodromy, identifying
$M_S$ with $h_S$.
\end{proof}

\subsection{Strata at infinity}
Assume now that $S\subset\Sigma_b^\infty$.  Let
\[
 i_E^\infty:Y^\infty:=\oY\setminus Y\hookrightarrow\oY
\]
denote the  inclusion.
As already mentioned, at a general point $x\in S$, the support of the sliced vanishing cycles 
$\varphi_{\og-b}^H(j_{E*}\Q_Y^H[d_S])$ is isolated.  Using the local product structure of $Y^\infty$ over the base, that is,
\[
 (Y^\infty,\og)\simeq
 \left((Y^\infty\cap\oG_b)\times D_b,\operatorname{pr}_2\right),
\]
we get that
\begin{equation}\label{eq28}
 \varphi_{\og-b}^H\left((i_E^\infty)_*
 \Q_{Y^\infty}^H[d_S]\right)=0
\end{equation}
near $x$.  Applying $\varphi_{\og-b}^H$ to the triangle
\[
 j_{E!}\Q_Y^H[d_S]\longrightarrow
 \Q_{\oY}^H[d_S]\longrightarrow
 (i_E^\infty)_*\Q_{Y^\infty}^H[d_S]
 \overset{[1]}{\lra}
\]
therefore yields the mixed Hodge module lift of the 
identification in \cite[Lemma~2.5]{MM26}:
\begin{equation}\label{eq29}
 \varphi_{\og-b}^H(j_{E!}\Q_Y^H[d_S])
 \simeq
 \varphi_{\og-b}^H(\Q_{\oY}^H[d_S])
\end{equation}
near $x$. We can now show the following.

\begin{theorem}
\label{thm2}
Let $S\subset\Sigma_b^\infty$ be a connected $s_b$-dimensional
stratum and let $x\in S$.  Then
\begin{align}
 (\Nu_S)_1
 &\simeq(H_S^\pitchfork)_1^\vee(-d_S),
 \label{eq30}\\
 (\Nu_S)_\lambda
 &\simeq(H_S^\pitchfork)_{\lambda^{-1}}^\vee(1-d_S),
 \quad\lambda\ne1,
 \label{eq31}
\end{align}
and these are isomorphisms of mixed Hodge structures with nilpotent
monodromy.  
\end{theorem}

\begin{proof} The proof follows the lines of \cite[Lemma 2.5]{MM26}, while keeping track of the Tate twists introduced by duality. 
 
Since $\cP_b^H|_S\simeq\cV_S[s_b]$, we have
\begin{equation}\label{eq32}
 \begin{aligned}
 \Nu_S
 =(\cV_S)_x
 =H^{-s_b}(i_x^*\cP_b^H)
 =H^0i_x^*\left(\oi_b^*\cP_b^H[-s_b]\right)
 \overset{\eqref{eq25}}{\simeq}
 H^0i_x^*\left(\varphi_{\og-b}^H
 \left(j_{E*}\Q_Y^H[d_S]\right)\right).
 \end{aligned}
\end{equation}

Since $Y$ is smooth of dimension $d_S$, self-duality implies
\[
 \D\left(\Q_Y^H[d_S]\right)
 \simeq\Q_Y^H[d_S](d_S),
\]
and therefore
\begin{equation}
 \D\left(j_{E!}\Q_Y^H[d_S]\right)
 \simeq j_{E*}\D\left(\Q_Y^H[d_S]\right)
 \simeq j_{E*}\Q_Y^H[d_S](d_S).
 \label{eq33}
\end{equation}
Equivalently,
\begin{equation}\label{eq34}
 j_{E*}\Q_Y^H[d_S]
 \simeq
 \D\left(j_{E!}\Q_Y^H[d_S]\right)(-d_S).
\end{equation}

Since the support of
$\varphi_{\og-b}^H(j_{E!}\Q_Y^H[d_S])$ is $\{x\}$ in a neighborhood
of $x$, its stalk and costalk at $x$ coincide, that is,
\begin{equation}\label{eq35}
 i_x^!\left(\varphi_{\og-b}^H(j_{E!}\Q_Y^H[d_S])\right)
 \xrightarrow{\sim}
 i_x^*\left(\varphi_{\og-b}^H(j_{E!}\Q_Y^H[d_S])\right),
\end{equation}
and both sides have cohomology only in degree zero.  Therefore,
\begin{align}
 H^0i_x^!\left(\varphi_{\og-b}^H(j_{E!}\Q_Y^H[d_S])\right)
 &\overset{\eqref{eq35}}{\simeq}
 H^0i_x^*\left(\varphi_{\og-b}^H(j_{E!}\Q_Y^H[d_S])\right)\notag\\
 &\overset{\eqref{eq29}}{\simeq}
 H^0i_x^*\left(\varphi_{\og-b}^H\Q_{\oY}^H[d_S]\right)\notag\\
 &\simeq
 \widetilde H^{d_S-1}(MF_{\og-b,x};\Q)\notag\\
 &=\widetilde H^{n-s_b}(MF_S^\pitchfork;\Q)=H_S^\pitchfork,
 \label{eq36}
\end{align}
where $MF_S^\pitchfork:=MF_{\og-b,x}$
denotes as before the transversal Milnor fiber to the stratum $S$. 
The isomorphism \eqref{eq36} respects the decomposition by
vanishing cycle monodromy, hence
\begin{align}
 H^0i_x^!\left(\varphi_{\og-b,1}^H(j_{E!}\Q_Y^H[d_S])\right)
 &\simeq(H_S^\pitchfork)_1,
 \label{eq37}\\
 H^0i_x^!\left(\psi_{\og-b,\mu}^H(j_{E!}\Q_Y^H[d_S])\right)
 &\simeq(H_S^\pitchfork)_\mu,\quad\mu\ne1.
 \label{eq38}
\end{align}

For the eigenvalue-one summand, we have
\begin{align*}
 (\Nu_S)_1
 &\overset{\eqref{eq32}}{\simeq}
 H^0 i_x^*\left(\varphi_{\og-b,1}^H
 \left(j_{E*}\Q_Y^H[d_S]\right)\right)\\
 &\overset{\eqref{eq34}}{\simeq}
 H^0 i_x^*\left(\varphi_{\og-b,1}^H
 \left(\D(j_{E!}\Q_Y^H[d_S])\right)\right)(-d_S)\\
 &\overset{\eqref{eq12}}{\simeq}
 H^0 i_x^*\left(\D\varphi_{\og-b,1}^H
 \left(j_{E!}\Q_Y^H[d_S]\right)\right)(-d_S)\\
 &\simeq
 \left(H^0 i_x^!\left(\varphi_{\og-b,1}^H
 \left(j_{E!}\Q_Y^H[d_S]\right)\right)\right)^\vee(-d_S)\\
 &\overset{\eqref{eq37}}{\simeq}
 (H_S^\pitchfork)_1^\vee(-d_S),
\end{align*}
which proves \eqref{eq30}.

For $\lambda\ne1$, the $\lambda$-summand of the 
vanishing cycle object is $\psi_{\og-b,\lambda}^H$.  Using
\eqref{eq13} and then
\eqref{eq38} with $\mu=\lambda^{-1}$, we obtain
\begin{align*}
 (\Nu_S)_\lambda
 &\overset{\eqref{eq32}}{\simeq}
 H^0i_x^*\left(\psi_{\og-b,\lambda}^H
 \left(j_{E*}\Q_Y^H[d_S]\right)\right)\\
 &\overset{\eqref{eq34}}{\simeq}
 H^0i_x^*\left(\psi_{\og-b,\lambda}^H
 \left(\D(j_{E!}\Q_Y^H[d_S])\right)\right)(-d_S)\\
 &\overset{\eqref{eq13}}{\simeq}
 H^0 i_x^*\left(\D
 \psi_{\og-b,\lambda^{-1}}^H
 \left(j_{E!}\Q_Y^H[d_S]\right)\right)(1-d_S)\\
 &\simeq
 \left(
 H^0i_x^!\left(\psi_{\og-b,\lambda^{-1}}^H
 \left(j_{E!}\Q_Y^H[d_S]\right)\right)
 \right)^\vee(1-d_S)\\
 &\overset{\eqref{eq38}}{\simeq}
 (H_S^\pitchfork)_{\lambda^{-1}}^\vee(1-d_S),
\end{align*}
which shows \eqref{eq31}.
\end{proof}

Combining Theorem~\ref{thm1},
Proposition~\ref{prop2}, and
Theorem~\ref{thm2}, we obtain  the following explicit form of Corollary \ref{cor1}.

\begin{corollary}
\label{cor2}
With $d_S=n-s_b+1$ and $\lambda=1$, we have:
\begin{equation} \label{eq39}
 h_1^{p,q}\left(V_H^{n-s_b}(b)\right)
 \le{}
 \sum_{\substack{S\subset\Sigma_b^a\\\dim S=s_b}}
 h_1^{p,q}(H_S^\pitchfork) 
 +
 \sum_{\substack{S\subset\Sigma_b^\infty\\\dim S=s_b}}
 h_1^{d_S-p,d_S-q}(H_S^\pitchfork).
\end{equation}
For $\lambda\ne1$,
\begin{equation} \label{eq40}
 h_\lambda^{p,q}\left(V_H^{n-s_b}(b)\right)
 \le{}
 \sum_{\substack{S\subset\Sigma_b^a\\\dim S=s_b}}
 h_\lambda^{p,q}(H_S^\pitchfork)
 +
 \sum_{\substack{S\subset\Sigma_b^\infty\\\dim S=s_b}}
 h_{\lambda^{-1}}^{d_S-1-p,d_S-1-q}(H_S^\pitchfork).
\end{equation}
Each term on the right can be sharpened by replacing it with the
local system fixed part of the corresponding variation.
\end{corollary}

\begin{proof}
For $S\subset\Sigma_b^a$, we get by evaluation and
\eqref{eq27} that
\[
 H^0(S;\cV_{S,\lambda})
 \hookrightarrow(\Nu_S)_\lambda
 \xrightarrow{\sim}(H_S^\pitchfork)_\lambda.
\]
For $S\subset\Sigma_b^\infty$, evaluation and
\eqref{eq30}--\eqref{eq31}
imply
\begin{align*}
 H^0(S;\cV_{S,1})
 &\hookrightarrow
 (H_S^\pitchfork)_1^\vee(-d_S),\\
 H^0(S;\cV_{S,\lambda})
 &\hookrightarrow
 (H_S^\pitchfork)_{\lambda^{-1}}^\vee(1-d_S),
 \quad\lambda\ne1.
\end{align*}
Using \eqref{eq6}, we get
\begin{align*}
 h^{p,q}\left((H_S^\pitchfork)_1^\vee(-d_S)\right)
 &=h_1^{d_S-p,d_S-q}(H_S^\pitchfork),\\
 h^{p,q}\left((H_S^\pitchfork)_{\lambda^{-1}}^\vee(1-d_S)\right)
 &=h_{\lambda^{-1}}^{d_S-1-p,d_S-1-q}(H_S^\pitchfork).
\end{align*}
Substituting into \eqref{eq24} proves
\eqref{eq39} and
\eqref{eq40}.
\end{proof}


\section{Limit mixed Hodge structure and monodromy on the cohomology of the general fiber}

Let $U=\C\setminus B$, choose $c\in U$, and set $F=f^{-1}(c)$ for the general fiber of $f$.  For
$b\in B$, we choose a positively oriented simple loop around $b$ and let
$T_b$ be the induced monodromy on $H^k(F;\Q)$.  The local system $(R^k f_*\Q_X)|_U$, whose stalk at $t\in U$
is $H^k(F_t;\Q)$, underlies an admissible variation of mixed
Hodge structures on $U$. For each $b\in B$ and $t\in U$ sufficiently close to $b$, the stalk
$\left(\psi_{t-b}(R^k f_*\Q_X)\right)_b$ is isomorphic to
$H^k(F_t;\Q)$ and
carries the limit mixed Hodge structure, which does not in general coincide
with Deligne's mixed Hodge structure on $H^k(F_t;\Q)$.

\begin{definition}\label{def2}
The \emph{limit mixed Hodge structure at $b$} is
\begin{equation}\label{eq41}
 \Hlim^k(F;b):=
 H^{k-n}\!\left(
 (\oa_b)_*\psi_{\of-b}^H\cM\right),
\end{equation}
where $\oa_b:\oF_b\to\mathrm{pt}$ is the constant map, $\psi^H$ denotes the total nearby cycle functor for mixed Hodge modules, and  $\cM:=j_*\Q_X^H[n+1]\in\MHM(\oX)$.
\end{definition}

We next recall from \cite{MM26} the identification of the underlying rational vector
space of \eqref{eq41} with $H^k(F;\Q)$.
Let $j:X\hookrightarrow\oX$ be the open inclusion, and denote by
$z$ the coordinate on $\C$. Since $f=\of\circ j$, we have
\[
 Rf_*\Q_X\simeq R\of_*Rj_*\Q_X.
\]
Applying to $\of$ the fact that proper direct image commutes with nearby cycles yields
\[
 \psi_{z-b}(Rf_*\Q_X)
 \simeq
 R(\of|_{\oF_b})_*
 \psi_{\of-b}(Rj_*\Q_X),
\]
where $\psi$ denotes the standard (unshifted) nearby cycle functor.  Thus
\[
 \bH^k\!\left(
   \oF_b;\psi_{\of-b}(Rj_*\Q_X)
 \right)
 \simeq
 H^k\!\left((\psi_{z-b}Rf_*\Q_X)_b\right).
\]
Furthermore, if $D_b$ is a small disk at $b$ such that $D_b\cap B=\{b\}$, and $D_b^*=D_b\setminus\{b\}$, 
then the cohomology sheaves of $(Rf_*\Q_X)|_{D_b^*}$ are local systems, with stalk at  $t\in D_b^*$ given by
\[
 H^k\!\left((\psi_{z-b}Rf_*\Q_X)_b\right)
 \simeq (R^k f_*\Q_X)_t.
\]

To identify this stalk with the cohomology of $F_t$, consider
the inclusions 
\[
 \oi_t:\oF_t\hookrightarrow\oX,
 \quad j_t:F_t\hookrightarrow\oF_t.
\]
After shrinking $D$, generic base change for $Rj_*$ yields for $t\in D_b^*$ that
\[
 \oi_t^*Rj_*\Q_X\simeq Rj_{t*}\Q_{F_t}.
\]
So we get by proper base change for $\of$ that
\[
 (Rf_*\Q_X)_t
 \simeq R\Gamma(\oF_t;\oi_t^*Rj_*\Q_X)
 \simeq R\Gamma(\oF_t;Rj_{t*}\Q_{F_t})
 \simeq R\Gamma(F_t;\Q).
\]
Combining these isomorphisms, we obtain
\[
 \bH^k\!\left(
   \oF_b;\psi_{\of-b}(Rj_*\Q_X)
 \right)
 \simeq H^k(F_t;\Q).
\]
Finally, a chosen path in $U$ from $t$ to $c$ identifies
$H^k(F_t;\Q)$ with $H^k(F;\Q)$ via the local system
$(R^k f_*\Q_X)|_U$. This then identifies the underlying rational
vector space of \eqref{eq41} with $H^k(F;\Q)$.
The mixed Hodge structure in \eqref{eq41} is the
limit mixed Hodge structure at $b$, but it does not in general agree,
under this identification, with Deligne's mixed Hodge
structure on $H^k(F;\Q)$.

  The semisimple
part of $T_b$ acts on the limit mixed Hodge structure, and on the
unipotent part
\[
 N_b:\Hlim^k(F;b)_1\longrightarrow
 \Hlim^k(F;b)_1(-1)
\]
is a morphism of mixed Hodge structures.

\smallskip

With the above notation, we can now show the following.
\begin{proposition}
\label{prop3}
For every $b\in B$ and $k\in\bZ$, we have:
\begin{enumerate}[label=\textup{(\roman*)},leftmargin=2.2em]
\item for $\lambda\ne1$, the canonical morphism induces an
isomorphism
\begin{equation}\label{eq42}
 \Hlim^k(F;b)_\lambda\xrightarrow{\sim}V_H^k(b)_\lambda.
\end{equation}
\item $\im N_b$ is a mixed Hodge subquotient of $V_H^k(b)_1$.
\item
\begin{equation}\label{eq43}
 \rk(T_b-\id\mid H^k(F;\Q))\le\dim V^k(b).
\end{equation}
\item if $k<n-s_b$, then $T_b=\id$ on
$\widetilde H^k(F;\Q)$.
\end{enumerate}
Moreover,
\begin{align}
 h_\lambda^{p,q}\left(\Hlim^k(F;b)\right)
 &=h_\lambda^{p,q}\left(V_H^k(b)\right), \ \ \lambda\ne1,
 \label{eq44}\\
 h^{p,q}(\im N_b)
 &\le h_1^{p,q}\left(V_H^k(b)\right).
 \label{eq45}
\end{align}
\end{proposition}

\begin{proof}
Let $K=Rf_*\Q_X$, and let
$i_b:\{b\}\hookrightarrow\C$ be the inclusion. The specialization
triangle
\[
 i_b^*K \longrightarrow
 \psi_{z-b}K
 \overset{\operatorname{can}}{\lra}
 \varphi_{z-b}K
 \overset{+1}{\lra}
\]
is compatible with monodromy, where the monodromy on $i_b^*K$
is the identity. Here both nearby and vanishing cycle functors are unshifted.
This triangle lifts to the derived category of mixed Hodge
modules and yields the exact sequence of mixed Hodge structures
\[
 H^k(i_b^*K)
 \longrightarrow \Hlim^k(F;b)
 \xrightarrow{\operatorname{can}} V_H^k(b)
 \longrightarrow H^{k+1}(i_b^*K).
\]
The two outer terms have only the eigenvalue $1$.
On the other hand, after complexification, taking generalized $\lambda$-eigenspaces
is an exact functor. Thus, for $\lambda\ne1$, the above sequence becomes
\[
 0\longrightarrow \Hlim^k(F;b)_\lambda
 \xrightarrow{\operatorname{can}_\lambda}
 V_H^k(b)_\lambda
 \longrightarrow 0.
\]
In particular, $\operatorname{can}_\lambda$ is an isomorphism
compatible with the induced Hodge and weight filtrations.
Over $\Q$, we get an isomorphism in $\MHM(\mathrm{pt})$:
\[
 \operatorname{can}_{\ne1}:
 \Hlim^k(F;b)_{\ne1}
 \xrightarrow{\sim} V_H^k(b)_{\ne1},
\]
where the subscript $\ne1$ denotes the non-unipotent summand, i.e., the rational summand complementary to the generalized eigenvalue-one summand.
In view of \eqref{eq42}, strictness of a morphism of mixed Hodge structures yields
\eqref{eq44}.

On the eigenvalue-one summand, the canonical and variation
morphisms are morphisms of mixed Hodge structures
\[
 \operatorname{can}:
 \Hlim^k(F;b)_1\longrightarrow V_H^k(b)_1,
 \quad
 \operatorname{Var}:
 V_H^k(b)_1\longrightarrow\Hlim^k(F;b)_1(-1),
\]
and satisfy
\[
 N_b=\operatorname{Var}\circ\operatorname{can}:
 \Hlim^k(F;b)_1\longrightarrow\Hlim^k(F;b)_1(-1).
\]
So $\operatorname{im}N_b$ gets a mixed
Hodge structure induced from that of the target
$\Hlim^k(F;b)_1(-1)$.

Let
\[
 I_b^k:=\operatorname{im}\!\left(
 \operatorname{can}:\Hlim^k(F;b)_1\longrightarrow V_H^k(b)_1
 \right).
\]
Since $\operatorname{can}$ is a morphism of mixed Hodge
structures, $I_b^k$ is a mixed Hodge substructure of
$V_H^k(b)_1$. Moreover,
\[
 \operatorname{Var}(I_b^k)
 =\operatorname{Var}\left(\operatorname{im}(\operatorname{can})\right)
 =\operatorname{im}(\operatorname{Var}\circ\operatorname{can})
 =\operatorname{im}N_b.
\]
Thus the restriction of $\operatorname{Var}$ to $I_b^k$ is
surjective onto $\operatorname{im}N_b$. Its kernel is
$I_b^k\cap\ker(\operatorname{Var})$, so we obtain a short exact
sequence of mixed Hodge structures
\[
 0\longrightarrow I_b^k\cap\ker(\operatorname{Var})
 \longrightarrow I_b^k
 \overset{\operatorname{Var}}{\lra}\operatorname{im}N_b
 \longrightarrow 0.
\]
In particular, the map
\[
 I_b^k/\left(I_b^k\cap\ker(\operatorname{Var})\right)
 \overset{\sim}{\lra}\operatorname{im}N_b,
 \quad [v]\mapsto\operatorname{Var}(v),
\]
is an isomorphism of mixed Hodge structures. This proves that
$\operatorname{im}N_b$ is a subquotient of $V_H^k(b)_1$.

To obtain the upper bounds on Hodge numbers, we complexify and apply
$\operatorname{Gr}_F^p\operatorname{Gr}_{p+q}^W$. Exactness of these functors on mixed Hodge
structures yields a surjection
\[
 \operatorname{Gr}_F^p\operatorname{Gr}_{p+q}^W(I_b^k)_\C
 \twoheadrightarrow
 \operatorname{Gr}_F^p\operatorname{Gr}_{p+q}^W
 (\operatorname{im}N_b)_\C
\]
and an injection
\[
 \operatorname{Gr}_F^p\operatorname{Gr}_{p+q}^W(I_b^k)_\C
 \hookrightarrow
 \operatorname{Gr}_F^p\operatorname{Gr}_{p+q}^W
 \left(V_H^k(b)_1\right)_\C.
\]
Taking dimensions yields
\[
 h^{p,q}(\operatorname{im}N_b)
 \leq h^{p,q}(I_b^k)
 \leq h_1^{p,q}\left(V_H^k(b)\right),
\]
which proves (ii) and \eqref{eq45}.

We next prove (iii). Denote as above $K=Rf_*\Q_X$ and
$i_b:\{b\}\hookrightarrow\C$ the inclusion.
Using the nearby fiber identification discussed before the statement of this proposition,
the specialization sequence becomes
\[
 H^k(i_b^*K)
 \xrightarrow{\operatorname{sp}}
 H^k(F;\Q)
 \xrightarrow{\operatorname{can}}
 V^k(b).
\]
The specialization map is compatible with monodromy, which
acts as the identity on $H^k(i_b^*K)$. Hence
$T_b\circ\operatorname{sp}=\operatorname{sp}$, and therefore
\[
 \ker(\operatorname{can})
 =\operatorname{im}(\operatorname{sp})
 \subseteq\ker(T_b-\operatorname{id}).
\]
By rank-nullity,
\begin{align*}
 \operatorname{rk}(T_b-\operatorname{id})
 &=\dim_\Q H^k(F;\Q)
   -\dim_\Q\ker(T_b-\operatorname{id})\\
 &\leq\dim_\Q H^k(F;\Q)
   -\dim_\Q\ker(\operatorname{can})\\
 &=\dim_\Q\operatorname{im}(\operatorname{can})\\
 &\leq\dim_\Q V^k(b),
\end{align*}
as desired.

We finally prove (iv). If $k<n-s_b$, then $V^k(b)=0$.
By (iii), we then get
\[
 0\leq\operatorname{rk}(T_b-\operatorname{id})
 \leq\dim_\Q V^k(b)=0.
\]
Hence $T_b=\operatorname{id}$ on $H^k(F;\Q)$, and therefore
also on $\widetilde H^k(F;\Q)$.
\end{proof}

In what follows, we let $T_b^{\mathrm{lim}}$ denote the nearby cycle monodromy on $\Hlim^k(F;b)$, while $M_b$ denotes as before  the vanishing cycle
monodromy on $V_H^k(b)$. Generalized eigenspaces will be taken
after complexification, with respect to the indicated monodromy.
The semisimple parts preserve the Hodge and weight filtrations,
so these filtrations restrict to the generalized eigenspaces.

For the eigenvalue-one summand, we use the morphism of mixed
Hodge structures
\[
 N_b:\Hlim^k(F;b)_1\longrightarrow\Hlim^k(F;b)_1(-1).
\]
In the statements below, $\operatorname{im}N_b$ always carries
the mixed Hodge structure induced from
$\Hlim^k(F;b)_1(-1)$. In the notation of Theorem \ref{thm1}, we then have the following.

\begin{corollary}
\label{cor3}
Assume $s_b>0$. For every $\lambda\ne1$,
\begin{equation}\label{eq46}
 h_\lambda^{p,q}\left(\Hlim^{n-s_b}(F;b)\right)
 \leq
 \sum_i h_\lambda^{p,q}
 \left(H^0(S_{b,i};\cV_{b,i})\right).
\end{equation}
On the eigenvalue-one summand,
\begin{equation}\label{eq47}
 h^{p,q}(\operatorname{im}N_b)
 \leq
 \sum_i h_1^{p,q}
 \left(H^0(S_{b,i};\cV_{b,i})\right).
\end{equation}
\end{corollary}

\begin{proof}
By \eqref{eq24},
\[
 h_\lambda^{p,q}\left(V_H^{n-s_b}(b)\right)
 \leq
 \sum_i h_\lambda^{p,q}
 \left(H^0(S_{b,i};\cV_{b,i})\right)
\]
for every $\lambda$.
Furthermore, by 
Proposition~\ref{prop3}, for $\lambda\ne1$,
\[
 h_\lambda^{p,q}\left(\Hlim^{n-s_b}(F;b)\right)
 =
 h_\lambda^{p,q}\left(V_H^{n-s_b}(b)\right),
\]
whereas
\[
 h^{p,q}(\operatorname{im}N_b)
 \leq h_1^{p,q}\left(V_H^{n-s_b}(b)\right).
\]
Altogether, these prove \eqref{eq46} and
\eqref{eq47}, respectively.
\end{proof}

As in \cite{MM26}, let us set
\[ s:=\max_{b\in B}s_b>0,\]
with
\[B_s=\{b\in B\mid s_b=s\}.\]

\begin{corollary}
\label{cor4}
On $\widetilde H^{n-s}(F;\Q)$, the monodromy $T_b$ is the identity
for every $b\in B\setminus B_s$. For $b\in B_s$, we have
\begin{equation}\label{eq49}
 \operatorname{rk}\!\left(
 T_b-\operatorname{id}\mid\widetilde H^{n-s}(F;\Q)
 \right)
 \leq\dim_\Q V^{n-s}(b)
 \leq\sum_i\dim_\Q H^0(S_{b,i};\cV_{b,i})
 \leq\sum_i\mu_{b,s,i}^\pitchfork.
\end{equation}
\end{corollary}

\begin{proof}
For $b\in B\setminus B_s$, we have $n-s<n-s_b$, so
Proposition~\ref{prop3}(iv) proves the first assertion.

For $b\in B_s$, the first two inequalities in
\eqref{eq49} follow from
Proposition~\ref{prop3}(iii) and
Theorem~\ref{thm1}, respectively.
For each $S=S_{b,i}$, evaluation at $x\in S$ is injective, hence
\[
 H^0(S;\cV_{b,i})\hookrightarrow\Nu_S.
\]
Proposition~\ref{prop2} identifies $\Nu_S$
with $H_S^\pitchfork$ on affine strata, while
Theorem~\ref{thm2} identifies its eigenspaces
with Tate twists of the dual eigenspaces of $H_S^\pitchfork$
on strata at infinity. In both cases,
\[
 \dim_\Q\Nu_S=\dim_\Q H_S^\pitchfork
 =\mu_{b,s,i}^\pitchfork.
\]
Taking dimensions and summing over $i$ proves the last inequality.
\end{proof}

The vanishing cohomology dimension can be used to bound the size of eigenvalue-one Jordan blocks of local monodromy. More precisely, we get the following.
\begin{corollary}
\label{cor5}
Let $\delta_b^k=\dim_\Q V^k(b)$. The number of Jordan blocks
of $T_b$ with eigenvalue $1$ and size greater than one is
at most $\delta_b^k$. Every Jordan block with eigenvalue $1$
has size at most $\delta_b^k+1$.
If $s>0$, and $b\in B_s$, these conclusions hold
with $\delta_b^{n-s}$ replaced by either
\[
 \sum_i\dim_\Q H^0(S_{b,i};\cV_{b,i})
 \quad\text{or}\quad
 \sum_i\mu_{b,s,i}^\pitchfork.
\]
\end{corollary}

\begin{proof}
We work over $\C$, which does not change the rank of
$T_b-\operatorname{id}$. Let $J_r(1)$ be a Jordan block with eigenvalue $1$
and size $r$. 
Subtracting the identity $I_r$ leaves a matrix of rank $r-1$, with 
 first column consisting of zeros, and its remaining $r-1$ columns
 linearly independent.

Let $r_1,\ldots,r_\nu\geq2$ be the sizes of the blocks with
eigenvalue $1$ and size greater than one. The ranks of the
blocks add, so their contribution to the rank of
$T_b-\operatorname{id}$ is
$\sum_{j=1}^{\nu}(r_j-1)$. All other blocks contribute
nonnegative ranks. Proposition~\ref{prop3}(iii)
therefore implies
\[
 \sum_{j=1}^{\nu}(r_j-1)
 \leq\operatorname{rk}(T_b-\operatorname{id})
 \leq\delta_b^k.
\]
Since each $r_j-1\geq1$, we obtain
\[
 \nu\leq\sum_{j=1}^{\nu}(r_j-1)\leq\delta_b^k.
\]
For each $j$ we also have
\[
 r_j-1
 \leq\sum_{\ell=1}^{\nu}(r_\ell-1)
 \leq\delta_b^k,
\]
and hence $r_j\leq\delta_b^k+1$. The same bound also applies to blocks of size one since $\delta_b^k\geq0$.

In degree $n-s$ and for $b\in B_s$, the dimension estimates
in \eqref{eq49} bound $\delta_b^{n-s}$ by either of
the two sums mentioned in the statement, thus  proving the final assertion.
\end{proof}

\subsection{Monodromy at infinity}

Let us now choose a sufficiently large disk in $\C$ containing $B$,
and a base point $c$ on its boundary. Let $T_\infty$ denote
the monodromy of the counterclockwise boundary circle on
$\widetilde H^{n-s}(F;\Q)$, where $F=f^{-1}(c)$.

We choose based loops around the points of $B$, each oriented
counterclockwise, whose product represents the boundary circle.
Denote $B=\{b_1,\ldots,b_r\}$, ordering the loops so that
the induced monodromy operators satisfy
\begin{equation}\label{eq50}
 T_\infty=T_{b_r}\cdots T_{b_1}.
\end{equation}
Notice that the counterclockwise large circle in the $z$-plane is clockwise
around $w=0$ in the coordinate $w=1/z$. Monodromy along a
positively oriented loop in the $w$-coordinate is therefore
$T_\infty^{-1}$. The next result gives upper bounds for the rank of monodromy at infinity.

\begin{corollary}
\label{cor6}
With $s>0$ and $B_s=\{b\in B\mid s_b=s\}$ as above, we have
\begin{align}
 &\operatorname{rk}\!\left(
 T_\infty-\operatorname{id}\mid\widetilde H^{n-s}(F;\Q)
 \right)
 \leq\sum_{b\in B_s}\dim_\Q V^{n-s}(b)
 \leq
 \sum_{b\in B_s}\sum_i
 \dim_\Q H^0(S_{b,i};\cV_{b,i})
 \leq
 \sum_{b\in B_s}\sum_i\mu_{b,s,i}^\pitchfork.
 \label{eq51}
\end{align}
The same inequalities hold with $T_\infty$ replaced by
$T_\infty^{-1}$.
\end{corollary}

\begin{proof}
Let $A$ and $B$ be automorphisms of a finite-dimensional
vector space $E$. Then it is a well known fact in linear algebra that
\[  \operatorname{rk}(AB-\operatorname{id}) \leq \operatorname{rk}(A-\operatorname{id}) + \operatorname{rk}(B-\operatorname{id}).\]
Applying this inequality $r$ times to the product in
\eqref{eq50}, acting on
$E=\widetilde H^{n-s}(F;\Q)$, yields
\[
 \operatorname{rk}(T_\infty-\operatorname{id})
 \leq\sum_{b\in B}\operatorname{rk}(T_b-\operatorname{id}).
\]
By Corollary~\ref{cor4}, $T_b$ is the identity
on $E$ whenever $b\notin B_s$, so the corresponding terms in the above sum are zero.
For each $b\in B_s$, we have by the same corollary that
\[
 \operatorname{rk}(T_b-\operatorname{id})
 \leq\dim_\Q V^{n-s}(b)
 \leq\sum_i\dim_\Q H^0(S_{b,i};\cV_{b,i})
 \leq\sum_i\mu_{b,s,i}^\pitchfork.
\]
Summing these inequalities over $b\in B_s$ proves
\eqref{eq51}.

Finally,
\[
 T_\infty^{-1}-\operatorname{id}
 =-T_\infty^{-1}(T_\infty-\operatorname{id}).
\]
Since multiplication by the invertible operator $-T_\infty^{-1}$
preserves the dimension of the image, we get
\[
 \operatorname{rk}(T_\infty^{-1}-\operatorname{id})
 =\operatorname{rk}(T_\infty-\operatorname{id}),
\]
which proves the assertion for the inverse monodromy.
\end{proof}


\section{On the Hodge theory of atypical fibers}

Recall that for $b \in B$, we denote $s_b=\dim \Sigma_b$ and $s=\max_{b\in B} s_b$. We have by \cite{Br03} (see also \cite[Theorem 4.1]{MM26}) that
\begin{equation}\label{eq52}
 \widetilde H^k(F_b;\Q)=0, \quad k<n-s-1.
\end{equation}
Furthermore, in \cite[Theorem 4.8]{MM26}, we also bound the Betti number of $F_b$ in degree $n-s-1$ by the sum of transversal Milnor numbers corresponding to the top-dimensional strata in
$\Sigma_b^\infty$. We now carry out that argument using mixed
Hodge modules, keeping track of the Tate twist in the
variation morphism.

All reduced cohomology groups below carry their usual mixed
Hodge structures. Our assumption $\dim\operatorname{Sing}(f)<n$
implies $s\leq n-1$, hence $n-s\geq1$.

\subsection{The unipotent morphism at infinity}

Let
\[
 i_b:F_b^\infty\hookrightarrow\oF_b,
 \quad
 j_b:F_b\hookrightarrow\oF_b
\]
be the complementary inclusions. On the fiber $\oF_b=\of^{-1}(b)$, consider the complex
\[
 \cF_b^H=\ok_b^*\cM[-n-1]\in D^b\MHM(\oF_b),
\]
where $\ok_b:\oF_b\hookrightarrow\oX$ is the inclusion and $\cM=j_*\Q_X^H[n+1]$. Thus
\[
 \rat(\cF_b^H)
 =\ok_b^*Rj_*\Q_X,
\]
and the underlying rational complexes of $\cF_b^H$ and
$\cP_b^H$ are the complexes denoted $\cF$ and $\cG$
in \cite[Section~4]{MM26}.
As before, we set
\[
 \cP_{b,1}^H=\varphi_{\of-b,1}^H\cM.
\]
This fits in the canonical and variation triangles in $D^b\MHM(\oF_b)$, that is,
\begin{equation}\label{eq53}
 \psi_{\of-b,1}^H\cM
 \overset{\operatorname{can}}{\lra}
 \cP_{b,1}^H
 \longrightarrow \ok_b^*\cM
 \overset{[1]}{\lra}
\end{equation}
and
\begin{equation}\label{eq54}
 \ok_b^!\cM
 \longrightarrow\cP_{b,1}^H
 \overset{\operatorname{Var}}{\lra}
 \psi_{\of-b,1}^H\cM(-1)
 \overset{[1]}{\lra}.
\end{equation}

Let $i:X^\infty\hookrightarrow\oX$ and
$k_b^\infty:F_b^\infty\hookrightarrow\oX_\infty$
be the inclusions. Since
$\ok_b\circ i_b=i\circ k_b^\infty$, we get
\[
 i_b^!\ok_b^!\cM
 \simeq(k_b^\infty)^!i^!\cM.
\]
Moreover, the identity $i^!j_*=0$ implies $i^!\cM=0$, and thus
\[
 i_b^!\ok_b^!\cM=0.
\]
Applying $i_b^!$ to the variation triangle \eqref{eq54} therefore yields
an isomorphism
\[
 i_b^!\operatorname{Var}:
 i_b^!\cP_{b,1}^H
 \overset{\sim}{\lra}
 i_b^!\psi_{\of-b,1}^H\cM(-1).
\]
After twisting by $(1)$, this becomes
\begin{equation}\label{eq55}
 (i_b^!\operatorname{Var})(1):
 i_b^!\cP_{b,1}^H(1)
 \overset{\sim}{\lra}
 i_b^!\psi_{\of-b,1}^H\cM.
\end{equation}

On vanishing cycles, the canonical and variation morphisms
satisfy
\[
 \operatorname{can}(-1)\circ\operatorname{Var}
 =N:
 \cP_{b,1}^H\longrightarrow\cP_{b,1}^H(-1).
\]
Twisting this identity by $(1)$ yields
\[
 \operatorname{can}\circ\operatorname{Var}(1)=N(1).
\]
We may therefore use the variation isomorphism \eqref{eq55} 
to replace the first term in the canonical triangle \eqref{eq53}.
Since $\ok_b^*\cM=\cF_b^H[n+1]$, the resulting triangle is
\begin{equation}\label{eq56}
 i_b^!\cP_{b,1}^H(1)
 \overset{N(1)}{\lra}
 i_b^!\cP_{b,1}^H
 \longrightarrow i_b^!\cF_b^H[n+1]
 \overset{[1]}{\lra}.
\end{equation}
In particular, the morphism replacing the sheaf-theoretic
map $M_b-\operatorname{id}$ in the proof of \cite[Theorem 4.1]{MM26} is
\begin{equation}\label{eq57}
 N(1):
 i_b^!\cP_{b,1}^H(1)
 \longrightarrow i_b^!\cP_{b,1}^H.
\end{equation}
The canonical morphism is an isomorphism on the
non-unipotent nearby and vanishing cycle summands, so those
summands do not contribute any term to $i_b^!\cF_b^H[n+1]$.

The following is the Hodge-theoretic version of \cite[Proposition 4.3]{MM26}.
\begin{lemma}\label{lem2}
There is a monomorphism of mixed Hodge structures
\begin{equation}\label{eq58}
 \widetilde H^{n-s-1}(F_b;\Q)
 \hookrightarrow
 \ker\!\left(
 N(1):\bH^{-s}(F_b^\infty;i_b^!\cP_{b,1}^H)(1)
 \longrightarrow
 \bH^{-s}(F_b^\infty;i_b^!\cP_{b,1}^H)
 \right).
\end{equation}
\end{lemma}

\begin{proof}
Let us start by recalling the following commutative diagram of inclusions:
\[
\xymatrix{
F_b \ar[d]_{k_b}\ar[r]^{j_b}& \oF_b \ar[d]^{\ok_b} & \ar[d]^{k^\infty_b}\ar[l]_{i_b} F_b^\infty \\
X  \ar[r]_{j}& \oX & \ar[l]^{i} X^\infty.
}
\]
The attaching triangle for $i_b$ and $j_b$ applied to $\cF_b^H=\ok_b^*j_*\Q^H_X$ is
\begin{equation}\label{eq59}
 (i_b)_!i_b^!\cF_b^H
 \longrightarrow\cF_b^H
 \longrightarrow(j_b)_*(j_b)^*\cF_b^H
 \overset{[1]}{\lra}.
\end{equation}
We have by the above commutative diagram that
\[
 (j_b)^*\cF_b^H\simeq k_b^*j^*j_*\Q_X^H \simeq k_b^*\Q_X^H\simeq\Q_{F_b}^H.
\]
Applying cohomology to \eqref{eq59} therefore yields an exact sequence of
mixed Hodge structures
\begin{equation}\label{eq60}
 \lra \bH^{n-s-1}(\oF_b;\cF_b^H)
 \longrightarrow H^{n-s-1}(F_b;\Q)
 \longrightarrow \bH^{n-s}(F_b^\infty;i_b^!\cF_b^H)\lra.
\end{equation}

On underlying rational vector spaces, proper base change
for $\of$ yields the identification
\[
 \bH^r(\oF_b;\rat(\cF_b^H))
 \simeq H^r((Rf_*\Q_X)_b)
 \simeq H^r(T(F_b);\Q),
\]
as in \cite[Section~4]{MM26}. Under this identification, the first
map in the exact sequence \eqref{eq60} is restriction from the tube to
the fiber. We will use this comparison to determine its
image and the vanishing of its source.

Under the above comparison, the morphism
\[
 R\Gamma(\oF_b;\rat(\cF_b^H))
 \longrightarrow R\Gamma(F_b;\Q)
\]
is identified with the restriction
\[
 R\Gamma(T(F_b);\Q)\longrightarrow R\Gamma(F_b;\Q).
\]
The attaching triangle therefore identifies
\[
 R\Gamma(F_b^\infty;
          \rat(i_b^!\cF_b^H))
 \simeq
 \operatorname{Cone}\!\left(
 R\Gamma(T(F_b);\Q)\longrightarrow R\Gamma(F_b;\Q)
 \right)[-1]
 \simeq R\Gamma(T(F_b),F_b;\Q).
\]
In particular, on underlying rational vector spaces, we have
\[
\bH^{n-s}(F_b^\infty;\rat(i_b^!\cF_b^H))
 \simeq H^{n-s}(T(F_b),F_b;\Q).
\]
Next, the reduced cohomology long exact sequence of the pair contains
\[
 \widetilde H^{n-s-1}(T(F_b);\Q)
 \longrightarrow \widetilde H^{n-s-1}(F_b;\Q)
 \longrightarrow H^{n-s}(T(F_b),F_b;\Q),
\]
and its first term vanishes by \cite[(4.1)]{MM26}. Thus the
connecting morphism from the attaching triangle induces
an injective morphism of mixed Hodge structures
\begin{equation}\label{eq61}
 \alpha_b:\widetilde H^{n-s-1}(F_b;\Q)
 \hookrightarrow \bH^{n-s}(F_b^\infty;i_b^!\cF_b^H).
\end{equation}
We also have by \cite[(4.10)]{MM26} the vanishing 
\[
 \bH^r(F_b^\infty;\rat(i_b^!\cP_b^H))=0,
 \quad r<-s.
\]
Since a mixed Hodge structure with zero underlying rational vector
space is zero, we obtain
\begin{equation}\label{eq62}
 \bH^r(F_b^\infty;i_b^!\cP_b^H)=0,
 \quad r<-s.
\end{equation}
So vanishing in the same range also holds for the eigenvalue-one direct summand.

The hypercohomology long exact sequence for the triangle 
\eqref{eq56} contains the terms
\begin{align*}
 \bH^{-s-1}(F_b^\infty;i_b^!\cP_{b,1}^H)
 \longrightarrow
 \bH^{-s-1}(F_b^\infty;i_b^!\cF_b^H[n+1])
 &\longrightarrow
 \bH^{-s}(F_b^\infty;i_b^!\cP_{b,1}^H)(1)\\
 &\overset{N(1)}{\lra}
 \bH^{-s}(F_b^\infty;i_b^!\cP_{b,1}^H),
\end{align*}
with first term vanishing by
\eqref{eq62}. Also,
\[
 \bH^{-s-1}(F_b^\infty;i_b^!\cF_b^H[n+1])
 =
 \bH^{n-s}(F_b^\infty;i_b^!\cF_b^H).
\]
We then get by exactness that
\[ \bH^{n-s}(F_b^\infty;i_b^!\cF_b^H) \cong 
 \ker\!\left(
 N(1):\bH^{-s}(F_b^\infty;i_b^!\cP_{b,1}^H)(1)
 \longrightarrow
 \bH^{-s}(F_b^\infty;i_b^!\cP_{b,1}^H)
 \right)
\]
as a mixed Hodge structure. Composing this isomorphism
with $\alpha_b$ proves \eqref{eq58}.
The resulting injection is a morphism of mixed Hodge
structures and is therefore strict.
\end{proof}

We next discuss Hodge bounds for an atypical fiber. 
Recall that $\cP_b^H=\varphi_{\of-b}^H\cM$, where
$\cM=j_*\Q_X^H[n+1]$ and $\varphi_{\of-b}^H$ is the 
vanishing cycle functor with underlying rational complex
$\varphi_{\of-b}[-1]$. For every connected $s$-dimensional
stratum $S\subset\Sigma_b^\infty$, the variation $\cV_S$
is defined by
\[
 \cP_b^H|_S\simeq\cV_S[s],
\]
and $\Nu_S=(\cV_S)_x$ denotes its stalk at a chosen point
$x\in S$. The subscript $1$ denotes the generalized
eigenvalue-one (unipotent) summand for vanishing cycle monodromy, so
$(\cV_{S,1})_x=(\Nu_S)_1$. 
The logarithm of the unipotent part of vanishing cycle
monodromy $M_b$ defines the morphism
\[
 N:\cP_{b,1}^H\longrightarrow\cP_{b,1}^H(-1)
\]
used above. Restricting this morphism to $S$ and shifting
by $[-s]$ yields
\[
 N:\cV_{S,1}\longrightarrow\cV_{S,1}(-1),
\]
which we denote by the same symbol.
Set
\begin{equation}\label{eq63}
 K_S:=H^0(S;\cV_{S,1}),
 \quad
 N_S:=H^0(S;N):
 K_S\rightarrow K_S(-1).
\end{equation}
Evaluation at $x$ identifies the underlying rational vector
space of $K_S$ with the subspace of $(\Nu_S)_1$ fixed by
the local system monodromy along $S$. Thus $N_S$ is induced
by logarithmic vanishing cycle monodromy, which should be distinguished from the local system (vertical) monodromy along the stratum $S$.

The following is a Hodge-theoretic enhancement of \cite[Corollary 4.7, Theorem 4.8]{MM26}.
\begin{theorem}
\label{thm3}
If $\dim\Sigma_b^\infty<s$, then
\[
 \widetilde H^{n-s-1}(F_b;\Q)=0.
\]
If $\dim\Sigma_b^\infty=s$, there is a monomorphism
of mixed Hodge structures
\begin{equation}\label{eq64}
 \widetilde H^{n-s-1}(F_b;\Q)
 \hookrightarrow
 \bigoplus_{\substack{S\subset\Sigma_b^\infty\\\dim S=s}}
 \ker\!\left(N_S(1):K_S(1)\longrightarrow K_S\right).
\end{equation}
Therefore,
\begin{align}
 h^{p,q}\left(\widetilde H^{n-s-1}(F_b;\Q)\right)
 &\leq
 \sum_S h^{p,q}
 \left(\ker(N_S(1):K_S(1)\to K_S)\right)
 \label{eq65}\\
 &\leq
 \sum_S h_1^{n-s-p,n-s-q}(H_S^\pitchfork).
 \label{eq66}
\end{align}
In particular,
\begin{equation}\label{eq67}
 \widetilde b_{n-s-1}(F_b)
 \leq\sum_S\dim_\Q\ker N_S
 \leq\sum_S\dim_\C(H_S^\pitchfork)_1
 \leq\sum_S\mu_S^\pitchfork.
\end{equation}
All sums are over the connected $s$-dimensional strata
contained in $\Sigma_b^\infty$.
\end{theorem}

\begin{proof}
Let $u_b^\infty:\Sigma_b^\infty\hookrightarrow F_b^\infty$
be the inclusion and consider 
\[
 \cK_b^H:=(u_b^\infty)^!i_b^!\cP_{b,1}^H
 \in D^b\MHM(\Sigma_b^\infty).
\]
Since $i_b^!\cP_{b,1}^H$ is supported on
$\Sigma_b^\infty$, the adjunction morphism yields an isomorphism
\[
 (u_b^\infty)_*\cK_b^H
 \xrightarrow{\sim}i_b^!\cP_{b,1}^H.
\]
Applying cohomology implies
\begin{equation}\label{eq68}
 \bH^{-s}(F_b^\infty;i_b^!\cP_{b,1}^H)
 \simeq \bH^{-s}(\Sigma_b^\infty;\cK_b^H).
\end{equation}

For $0\leq\ell\leq s$, let
\[
 U_{b,\infty}^\ell=
 \bigcup_{\substack{S\subset\Sigma_b^\infty\\\dim S\geq\ell}}S,
 \quad
 \Sigma_{b,\infty}^\ell=
 \bigcup_{\substack{S\subset\Sigma_b^\infty\\\dim S=\ell}}S,
\]
with inclusion maps
\[
 v_\ell^\infty:\Sigma_{b,\infty}^\ell
 \hookrightarrow U_{b,\infty}^\ell,
 \quad
 u_\ell^\infty:U_{b,\infty}^{\ell+1}
 \hookrightarrow U_{b,\infty}^\ell.
\]
Denote the restriction by
\[
 \cK_b^{H,\ell}
 :=\cK_b^H|_{U_{b,\infty}^\ell},
\]
and consider the attaching triangle 
\begin{equation}\label{eq69}
 (v_\ell^\infty)_!(v_\ell^\infty)^!
 \cK_b^{H,\ell}
 \longrightarrow\cK_b^{H,\ell}
 \longrightarrow
 (u_\ell^\infty)_*\cK_b^{H,\ell+1}
 \xrightarrow{[1]}.
\end{equation}

As in the proof of \cite[Theorem~4.6]{MM26}, the
costalk condition for the perverse sheaf $\rat(\cP_{b,1}^H)$
and the hypercohomology spectral sequence imply
\[
 \bH^r\!\left(\Sigma_{b,\infty}^\ell;
 (v_\ell^\infty)^!\cK_b^{H,\ell}\right)=0,
 \quad r<-\ell,
\]
where vanishing in mixed Hodge structures follows from
faithfulness of $\rat$.

Now, for $0\leq\ell<s$, the hypercohomology long exact sequence of the attaching
triangle \eqref{eq69} contains the terms
\begin{align*}
 \bH^{-s}\!\left(\Sigma_{b,\infty}^\ell;
 (v_\ell^\infty)^!\cK_b^{H,\ell}\right)
 \longrightarrow
 \bH^{-s}(U_{b,\infty}^\ell;\cK_b^{H,\ell}) \longrightarrow
 \bH^{-s}(U_{b,\infty}^{\ell+1};
          \cK_b^{H,\ell+1}),
\end{align*}
where the first term vanishes because $-s<-\ell$.
Thus each restriction map in this sequence is injective.
Composing them yields an $N$-equivariant monomorphism
\[
 \bH^{-s}(\Sigma_b^\infty;\cK_b^H)
 \hookrightarrow
 \bH^{-s}(U_{b,\infty}^s;\cK_b^{H,s}).
\]
When $s=0$, this map is the identity and no iteration is needed.
If $\dim\Sigma_b^\infty<s$, then
$U_{b,\infty}^s=\emptyset$, so the target is zero.
Together with \eqref{eq68} and Lemma~\ref{lem2}, this proves
the first assertion.

Suppose now that $\dim\Sigma_b^\infty=s$, hence $s_b=s$.
The strata in $U_{b,\infty}^s$ are therefore top-dimensional
strata also of $\Sigma_b$, so their union is open in $\Sigma_b$.
Regard $\cP_{b,1}^H$ as an object on its support $\Sigma_b$,
and let
\[
 \iota:U_{b,\infty}^s\hookrightarrow\Sigma_b
\]
be the inclusion. Since $\iota$ is open, $\iota^!=\iota^*$.
By the definition of $\cK_b^{H,s}$ and transitivity
of the functors $(-)^!$, we get
\[
 \cK_b^{H,s}
 \simeq\iota^!\cP_{b,1}^H
 \simeq\iota^*\cP_{b,1}^H
 =\cP_{b,1}^H|_{U_{b,\infty}^s}.
\]
On each connected
stratum $S\subset U_{b,\infty}^s$, we thus have
\[
 \cK_b^{H,s}|_S\simeq\cV_{S,1}[s].
\]
The strata $S$ are the open and closed components of
$U_{b,\infty}^s$, hence
\begin{align}
 \bH^{-s}(\Sigma_b^\infty;\cK_b^H)
 &\hookrightarrow
\bH^{-s}(U_{b,\infty}^s;\cK_b^{H,s})
 \simeq
 \bigoplus_{\substack{S\subset\Sigma_b^\infty\\\dim S=s}}
 H^0(S;\cV_{S,1})
 =\bigoplus_{\substack{S\subset\Sigma_b^\infty\\\dim S=s}}K_S.
 \label{eq70}
\end{align}

All restriction maps commute with $N$. Twisting the last
injection by $(1)$ therefore yields the commutative square
\[
 \begin{array}{ccc}
 \bH^{-s}(\Sigma_b^\infty;\cK_b^H)(1)
 &\lhook\joinrel\longrightarrow&\displaystyle\bigoplus_S K_S(1)\\[2mm]
 {\scriptstyle N(1)}\downarrow
 &&\downarrow{\scriptstyle\oplus_S N_S(1)}\\[1mm]
 \bH^{-s}(\Sigma_b^\infty;\cK_b^H)
 &\lhook\joinrel\longrightarrow&\displaystyle\bigoplus_S K_S.
 \end{array}
\]
If an element of the upper-left group belongs to
$\ker N(1)$, commutativity shows that its image in the
upper-right group belongs to $\ker(\oplus_S N_S(1))$.
The induced map on kernels is injective because the upper
horizontal map is injective. Since a direct-sum map acts
componentwise,
\[
 \ker(\oplus_S N_S(1))=\bigoplus_S\ker N_S(1).
\]
Combining this injection with
\eqref{eq68} and
Lemma~\ref{lem2} proves
\eqref{eq64}.
It is strict because it is a morphism of mixed Hodge structures.

For the local Hodge estimate, evaluation at $x\in S$ yields
\begin{equation}\label{eq71}
 \ker N_S(1)\hookrightarrow K_S(1)
 \hookrightarrow(\Nu_S)_1(1)
 \simeq(H_S^\pitchfork)_1^\vee(-n+s),
\end{equation}
where the last isomorphism follows from
Theorem~\ref{thm2} after twisting by $(1)$,
since $s_b=s$. Strictness then implies
\[
 h^{p,q}(\ker N_S(1))
 \leq h^{p,q}\left((H_S^\pitchfork)_1^\vee(-n+s)\right)
 =h_1^{n-s-p,n-s-q}(H_S^\pitchfork).
\]
Summing over the $s$-dimensional strata $S\subset\Sigma_b^\infty$ and using the strictness of \eqref{eq64},
proves \eqref{eq65} and
\eqref{eq66}.

Finally, duality and Tate twists preserve dimension, so the
inclusions in \eqref{eq71} imply that
\[
 \dim_\Q\ker N_S
 \leq\dim_\C(H_S^\pitchfork)_1
 \leq\dim_\Q H_S^\pitchfork
 =\mu_S^\pitchfork.
\]
Taking dimensions in \eqref{eq64}
and summing over the $s$-dimensional strata
$S\subset\Sigma_b^\infty$ proves \eqref{eq67}.
\end{proof}

\begin{remark}
The Milnor number bounds obtained in
\cite[Theorem~4.8]{MM26} are sharp. However, for a
particular polynomial, knowledge of the monodromy actions may
give smaller upper bounds. More precisely, for each
$s$-dimensional stratum $S\subset\Sigma_b^\infty$, we have
\[
 \dim_\Q\ker N_S
 \leq\dim_\Q K_S
 \leq\dim_\Q(\Nu_S)_1
 =\dim_\C(H_S^\pitchfork)_1
 \leq\mu_S^\pitchfork.
\]
After summing over these strata, the inequalities improve the
Milnor number bound whenever the eigenvalue-one condition,
invariance under monodromy along a stratum, or passage to the Hodge-theoretic kernel $\ker N_S$ decreases a summand. The Hodge number inequalities
also bound the individual Hodge types, which the
Milnor number bounds do not distinguish.
\end{remark}


\section{Hodge-filtered semicontinuity and applications}
Let $f:\C^{n+1}\to\C$ be a polynomial of degree $d\geq 2$, with $n\geq 1$,
and let $F$ be a general fiber. The top Betti number of $F$, $b_n(F)$,  is bounded above by 
$(d-1)^{n+1}$ and this bound is sharp, see \cite[Proposition 2.1]{ST07} and also \cite[Proposition~5.1]{MM26}. In this section, we obtain a refined sharp dimension upper bound for each graded piece of
the Hodge filtration on $H^n(F;\C)$. 
As in \cite{MM26}, this is a consequence of the semicontinuity property of the dimension of each such Hodge graded piece, which we prove in Theorem \ref{thm4} below, together with specific calculations in the case of a generic-at-infinity polynomial in the sense of \cite{T07}.

\medskip

Let $P(x,\tau)$ be a one-parameter polynomial family of
degree $d$ polynomials on $\C^{n+1}$, and let
$f_\tau=P(-,\tau)$. As in \cite[Theorem~5.2]{MM26}, choose
$c\in\C$ sufficiently general and a sufficiently small
parameter disk $\Delta$ centered at $0$ so that
\[
 F^\tau=f_\tau^{-1}(c)
\]
is a general fiber for every $\tau\in\Delta$.
In particular, each $F^\tau$ is smooth of dimension $n$.

Let $P^h$ be the degree $d$ homogenization of $P$ in the
$x$-variables, and define
\[
 \oZ=
 \left\{([x_0:x],\tau)\in\bP^{n+1}\times\C
 \ \middle|\
 P^h(x_0,x,\tau)-cx_0^d=0\right\}.
\]
The projection $\pi:\oZ\to\C_\tau$ to the $\tau$-factor is projective, and
for each $\tau\in\Delta$ its fiber
$\oF^\tau:=\pi^{-1}(\tau)$ is the projective closure
of $F^\tau$ in $\bP^{n+1}$. The affine part
\[
 Z:=\oZ\cap\{x_0\ne0\}
 \simeq\{(x,\tau)\in\C^{n+1}\times\C\mid P(x,\tau)=c\}
\]
then satisfies $F^\tau=\oF^\tau\cap Z$. Denote the inclusions by
\[
 u:Z\hookrightarrow\oZ,
 \quad
 \oi_0:\oF^0\hookrightarrow\oZ,
\]
and let $\pi_0:\oF^0\to\{0\}$ be the restriction of
$\pi$ to its central fiber.

Let us set
\[
 M_{\oZ}:=
 u_!\Q_Z^H[n+1].
\]
Since the hypersurface $Z$ is a local complete
intersection of pure dimension $n+1$, the complex 
$\Q_Z[n+1]$ is perverse. Therefore, 
$\Q_Z^H[n+1]\in \MHM(Z)$. The inclusion $u$ is an affine
open immersion and extension by zero along such an
inclusion is $t$-exact, hence
$M_{\oZ}\in\MHM(\oZ)$.

\begin{definition}\label{def3}
The compactly supported limit mixed Hodge structure is defined as
\[
 H_{c,\mathrm{lim}}^n(F^\tau;0):=
 H^0\!\left((\pi_0)_*\psi_\pi^H M_{\oZ}\right),
\]
where $\psi_\pi^H$ is the (total) nearby cycle functor on mixed Hodge modules.
\end{definition}

The above limit is taken in the parameter $\tau$ at $0$.
After choosing a nearby fiber identification, its underlying
rational vector space is $H_c^n(F^\tau;\Q)$ for
$0<|\tau|\ll1$, an identification which we will discuss in detail in the proof of the semicontinuity result below.

\begin{theorem}
\label{thm4}
There is a monomorphism of mixed Hodge structures
\begin{equation}\label{eq72}
 H_c^n(F^0;\Q)
 \hookrightarrow H_{c,\mathrm{lim}}^n(F^\tau;0).
\end{equation}
Therefore,
\begin{equation}\label{eq73}
 h^{p,q}\left(H_c^n(F^0;\Q)\right)
 \leq
 h^{p,q}\left(H_{c,\mathrm{lim}}^n(F^\tau;0)\right)
\end{equation}
for all $p,q$. Moreover, for every $p$ and every sufficiently small
$\tau\ne0$, we have
\begin{align}
 \dim_\C\operatorname{Gr}_F^pH_c^n(F^0;\C)
 &\leq
 \dim_\C\operatorname{Gr}_F^pH_c^n(F^\tau;\C),
 \label{eq74}\\
 \dim_\C\operatorname{Gr}_F^pH^n(F^0;\C)
 &\leq
 \dim_\C\operatorname{Gr}_F^pH^n(F^\tau;\C).
 \label{eq75}
\end{align}
In particular, $b_n(F^0)\leq b_n(F^\tau)$.
\end{theorem}

\begin{proof}
We start with the specialization triangle 
\begin{equation}\label{eq76}
 \oi_0^*M_{\oZ}[-1]
 \longrightarrow\psi_\pi^H M_{\oZ}
 \longrightarrow\varphi_\pi^H M_{\oZ}
 \overset{[1]}{\lra}
\end{equation}
and first prove the vanishing needed for injectivity.
Since $\pi$ is proper, we get by proper base change for vanishing cycles:
\[
 (\pi_0)_*\varphi_\pi^H M_{\oZ}
 \simeq
 \phi_{\id_\C}^H(\pi_*M_{\oZ}),
\]
where $\id_\C:\C_\tau\to\C$ is the coordinate function, whose
zero fiber is $\{0\}$. We also have by properness that
$\pi_*=\pi_!$, hence
\[
 \pi_*M_{\oZ}
 =\pi_*u_!\Q_Z^H[n+1]
 \simeq
 (\pi\circ u)_!\Q_Z^H[n+1].
\]
Therefore
\begin{equation}\label{eq77}
 (\pi_0)_*\varphi_\pi^H M_{\oZ}
 \simeq
 \varphi_{\id_\C}^H
 \left((\pi\circ u)_!\Q_Z^H[n+1]\right).
\end{equation}

The map $\pi\circ u:Z\to\C_\tau$ is affine, since 
$Z$ is a closed subvariety of
$\C^{n+1}\times\C_\tau$ and the projection to $\C_\tau$
is affine. It follows that 
\[
 (\pi\circ u)_!\Q_Z^H[n+1]
 \in D^{\geq0}\MHM(\C_\tau).
\]
Indeed, the corresponding statement for the underlying rational
complex is the left $t$-exactness of $R(\pi\circ u)_!$. 
Since $\rat$ is faithful
and exact on mixed Hodge modules, the vanishing of the
negative cohomology objects also holds in $D^b\MHM$.

The functor $\varphi_{\id_\C}^H$ is $t$-exact, so applying it to
the above object and using
\eqref{eq77} yields
\[
 (\pi_0)_*\varphi_\pi^H M_{\oZ}
 \in D^{\geq0}\MHM(\{0\}).
\]
In particular,
\[
 H^{-1}\!\left((\pi_0)_*\varphi_\pi^H M_{\oZ}\right)=0.
\]

Applying $(\pi_0)_*$ to
\eqref{eq76} and taking cohomology yields the exact sequence of mixed Hodge structures
\begin{align*}
 H^{-1}\!\left((\pi_0)_*\varphi_\pi^H M_{\oZ}\right)
 \longrightarrow
 H^0\!\left((\pi_0)_*\oi_0^*M_{\oZ}[-1]\right) \longrightarrow
 H^0\!\left((\pi_0)_*\psi_\pi^H M_{\oZ}\right).
\end{align*}
As the first group is zero, the second arrow is a
monomorphism of mixed Hodge structures:
\begin{equation}\label{eq78}
 H^0\!\left((\pi_0)_*\oi_0^*M_{\oZ}[-1]\right)
 \hookrightarrow
 H^0\!\left((\pi_0)_*\psi_\pi^H M_{\oZ}\right).
\end{equation}

We next identify the two mixed Hodge structures in \eqref{eq78}. Let
\[
 u_0:F^0\hookrightarrow\oF^0,
 \quad
 i_0:F^0\hookrightarrow Z
\]
be the inclusions. The square formed by
$u,u_0,\oi_0,i_0$ is Cartesian, so 
base change yields
\[
 \oi_0^*u_!
 \simeq(u_0)_!i_0^*.
\]
Since $i_0^*\Q_Z^H=\Q_{F^0}^H$, we then obtain
\[
 \oi_0^*M_{\oZ}[-1]
 \simeq(u_0)_!\Q_{F^0}^H[n].
\]
Since $\pi_0$ is proper,  $(\pi_0)_*=(\pi_0)_!$.
So letting $a_{F^0}:F^0\to\mathrm{pt}$ be the constant map, we  have
\begin{align*}
 H^0\!\left((\pi_0)_*\oi_0^*M_{\oZ}[-1]\right)
 &\simeq
 H^0\!\left((\pi_0)_!(u_0)_!\Q_{F^0}^H[n]\right)\simeq
 H^0\!\left((a_{F^0})_!\Q_{F^0}^H[n]\right)=H_c^n(F^0;\Q).
\end{align*}
These identifications are in the category of mixed Hodge
structures. On the other hand, the target of
\eqref{eq78} is
$H_{c,\mathrm{lim}}^n(F^\tau;0)$ by definition.
This proves \eqref{eq72}.

Let us now identify the rational vector space underlying the
target of 
\eqref{eq78} and explain its limit mixed Hodge structure.
Since proper direct image also commutes with nearby cycles, 
\[
 (\pi_0)_*\psi_\pi^H M_{\oZ}
 \simeq
 \psi_{\id_\C}^H
 \left((\pi\circ u)_!\Q_Z^H[n+1]\right).
\]
After applying $H^0$ and $\rat$, we get
\[
 \rat \ 
 H_{c,\mathrm{lim}}^n(F^\tau;0)
 \simeq
 H^n\!\left(
 \left(\psi_{\id_\C} R(\pi\circ u)_!\Q_Z\right)_0
 \right),
\]
where $\psi_{\id_\C}$ is the unshifted nearby cycle
functor.

Shrinking $\Delta$ if necessary, we can assume that the cohomology sheaves of
$R(\pi\circ u)_!\Q_Z$ are local systems
on $\Delta^*=\Delta\setminus\{0\}$.
We then have that
\[
 H^n\!\left(
 \left(\psi_{\id_\C} R(\pi\circ u)_!\Q_Z\right)_0
 \right)
 \simeq
 \left(R^n(\pi\circ u)_!\Q_Z\right)_\tau \simeq H_c^n(F^\tau;\Q),
\]
where the last identification follows by base change.
Thus the rational vector space underlying the defined limit
mixed Hodge structure is the  compactly supported
cohomology of a nearby fiber.
Moreover, the local system
\[
 \left(R^n(\pi\circ u)_!\Q_Z\right)|_{\Delta^*}
\]
underlies the admissible variation of mixed Hodge structure
induced by the mixed Hodge module functor $(\pi\circ u)_!$.
Its stalk at $\tau$ carries Deligne's mixed Hodge structure
on $H_c^n(F^\tau;\Q)$, and its limit at $0$ is precisely
$H_{c,\mathrm{lim}}^n(F^\tau;0)$ defined above.

Since \eqref{eq72} is a morphism of mixed Hodge
structures, it is strict for the Hodge and weight filtrations.
Applying
$\operatorname{Gr}_F^p\operatorname{Gr}_{p+q}^W$
after complexification therefore yields an injective map, and taking dimensions proves \eqref{eq73}.
Applying only $\operatorname{Gr}_F^p$ yields
\[
 \dim_\C\operatorname{Gr}_F^pH_c^n(F^0;\C)
 \leq
 \dim_\C\operatorname{Gr}_F^p
 H_{c,\mathrm{lim}}^n(F^\tau;0)_\C.
\]
To deduce \eqref{eq74}, we need to compare the
Hodge-graded dimensions of the limit
Hodge filtration with those of a nearby fiber.
We first make a finite ramified base change so that the monodromy is
unipotent. This does not change the nearby or limit
Hodge-graded dimensions. By admissibility, the Hodge
filtration extends by subbundles to Deligne's canonical
extension, with fiber at $0$ equal to the limit Hodge
filtration. The graded bundles have constant rank across the disk, so
\[
 \dim_\C\operatorname{Gr}_F^p
 H_{c,\mathrm{lim}}^n(F^\tau;0)_\C
 =
 \dim_\C\operatorname{Gr}_F^pH_c^n(F^\tau;\C).
\]
Combining this equality with the inequality obtained from
the strict monomorphism \eqref{eq72} proves
\eqref{eq74}.

Finally, each $F^\tau$, including $F^0$, is smooth of
complex dimension $n$. Since Poincar\'e duality is an isomorphism
of mixed Hodge structures, we get
\[
 H_c^n(F^\tau;\Q)
 \simeq H^n(F^\tau;\Q)^\vee(-n).
\]
Therefore, 
\[
 \operatorname{Gr}_F^p
 \left(H^n(F^\tau;\Q)^\vee(-n)\right)_\C
 \simeq
 \left(\operatorname{Gr}_F^{n-p}H^n(F^\tau;\C)\right)^\vee.
\]
Taking dimensions yields
\[
 \dim_\C\operatorname{Gr}_F^pH_c^n(F^\tau;\C)
 =
 \dim_\C\operatorname{Gr}_F^{n-p}H^n(F^\tau;\C).
\]
Applying this equality at $\tau=0$ and at $\tau\ne0$ in
\eqref{eq74}, and then replace $p$ by $n-p$, 
proves \eqref{eq75}.
The Hodge filtration is finite and exhaustive, so the sum
of the dimensions of its graded pieces is the dimension
of the cohomology group. Summing
\eqref{eq75} over $p$ therefore implies
$b_n(F^0)\leq b_n(F^\tau)$, a result already proved in \cite[Theorem 5.2]{MM26}.
\end{proof}

Recall that a polynomial is called {\it generic-at-infinity} if the projective
closure of each fiber is smooth near the hyperplane at infinity
$H_\infty\subset \bP^{n+1}$ and meets $H_\infty$ transversely, e.g., see \cite{T07,ST11}. The fibers of such polynomial can have at most isolated singularities. In particular, the cohomology of their general fiber is concentrated in middle degree $n$. 
The following observation reduces the calculation of Hodge numbers on the cohomology of the general fiber of a generic-at-infinity polynomial to the case of the Fermat polynomial
\[ P_F=x_0^d+\cdots +x_n^d\]
and its general fiber $F'=\{P_F=1\}$.

\begin{lemma}\label{lem4}
If $g:\C^{n+1}\to\C$ is a degree $d$ generic-at-infinity polynomial, then
the general fiber of $g$ has the same mixed Hodge numbers as $F'$.
\end{lemma}

\begin{proof}
Let $F$ be a general fiber of $g$, with projective closure $\oF \subset \bP^{n+1}$, and set $F^\infty=\oF\cap H_\infty$. The fiber $F$ is smooth and has the homotopy type of a bouquet of $n$-spheres, see 
\cite[Theorem~3.1]{ST95}. By the hypothesis at infinity, $\oF$ is also smooth along $F^\infty$, and $F^\infty$ is a smooth hypersurface of degree $d$ in $H_\infty \simeq \bP^n$. So $\oF$ is a smooth degree $d$ hypersurface in $\bP^{n+1}$.

Set $H^k_\infty(\oF)=H^k(\oF,F)$.
The cohomology long exact sequence of the pair $(\oF,F)$ is a sequence of mixed Hodge structures containing the following exact segment
\[
 0\longrightarrow H^n_\infty(\oF)\longrightarrow H^n(\oF)
 \longrightarrow H^n(F)\longrightarrow H^{n+1}_\infty(\oF)
 \longrightarrow H^{n+1}(\oF)\longrightarrow 0.
\]
Indeed, $H^{n+1}(F)=0$ since $F$ is affine of complex dimension $n$. At the other end, $H^{n-1}(F)=0$ for $n>1$ by the bouquet decomposition, while 
for $n=1$ the restriction $H^0(\oF)\to H^0(F)$ is an isomorphism
since both curves are connected.

By Lefschetz duality (see, e.g., \cite[Corollary~6.28]{PS08}), we have a mixed Hodge structure isomorphism
\[H^k_\infty(\oF)\simeq H^{2n-k}(F^\infty)^\vee(-n).\]
Since $\oF$ and $F^\infty$ are smooth projective hypersurfaces,
their Hodge numbers are known to depend only on their dimensions and on $d$. By additivity in the exact sequence, the same holds for the mixed Hodge numbers of $H^n(F)$. This applies in particular to $F'$, since $P_F$ is generic-at-infinity. \end{proof}

Let us now recall from \cite{Ste77} the computation of mixed Hodge numbers for the general fiber $F'=\{P_F=1\}$ of the Fermat polynomial $P_F$.
Consider the graded Milnor algebra of $P_F$,
\[
 M(P_F)=\C[x_0,\ldots,x_n]/(x_0^{d-1},\ldots,x_n^{d-1}).
\]
A basis of $M(P_F)$ consists of the monomials $x^\alpha=x_0^{\alpha_0}\cdots x_n^{\alpha_n}$
with $0\leq \alpha_j\leq d-2$. To each $x^\alpha$, 
Steenbrink associates the number
\[
 \ell(\alpha)=\frac{(\alpha_0+1)+\cdots+(\alpha_n+1)}{d}.
\]
By \cite[Theorem~1]{Ste77}, we know that the only possible weights on
$H^n(F')$ are $n$ and $n+1$. Moreover, for fixed $p$, the monomials with
$p<\ell(\alpha)<p+1$ index a basis of $\Gr_F^p\Gr_n^W H^n(F';\C)$, while those 
with $\ell(\alpha)=p$ index a basis of $\Gr_F^p\Gr_{n+1}^W H^n(F';\C)$.
In terms of the grading of the Milnor algebra, this means that
\begin{align*}
 h^{p,n-p}(H^n(F'))   &=\sum_{k=1}^{d-1}\dim_{\C} M(P_F)_{pd-(n+1)+k},\\
 h^{p,n+1-p}(H^n(F')) &=\dim_{\C} M(P_F)_{pd-(n+1)}.
\end{align*}
Adding the two rows yields
\begin{equation}\label{eq81}
 \dim_{\C}\Gr_F^p H^n(F';\C)
 =\sum_{k=0}^{d-1}\dim_{\C} M(P_F)_{pd-(n+1)+k}.
\end{equation}

Together with the semicontinuity of Theorem \ref{thm4}, we can now deduce the following sharp upper bound for the graded Hodge pieces of top cohomology of the general fiber for any polynomial. The bound is attained by the Fermat polynomial itself, or any generic-at-infinity polynomial (cf. Lemma \ref{lem4}).

\begin{proposition}\label{prop4}
Let $f:\C^{n+1}\to\C$ be a polynomial of degree $d\geq 2$, with $n\geq 1$,
and let $F$ be a general fiber. For $0 \leq p \leq n$, we have
\begin{equation}\label{eq82}
 \dim_{\C}\Gr_F^p H^n(F;\C) \leq \sum_{k=0}^{d-1}\dim_{\C} M(P_F)_{pd-(n+1)+k},
 \end{equation} 
 with $M(P_F)$ denoting as above the graded Milnor algebra of the degree $d$ Fermat polynomial $P_F=x_0^d+\cdots +x_n^d$.
\end{proposition}

\begin{proof}
Let $f_\tau=f+\tau P_F$. By \cite[Proposition~2.2]{ST11},
$f_\tau$ is generic-at-infinity for $0<|\tau|\ll 1$. Choose $c$ sufficiently
general so that $F^\tau=f_\tau^{-1}(c)$ is a general fiber for all
$\tau$ in a small disk about $0$, and let $F=F^0$ be the general fiber of $f=f_0$. By the Hodge-filtered semicontinuity \eqref{eq75} of Theorem \ref{thm4} together with Lemma~\ref{lem4}, we get
\[
 \dim_{\C}\Gr_F^p H^n(F;\C)
 \leq \dim_{\C}\Gr_F^p H^n(F^\tau;\C)
 = \dim_{\C}\Gr_F^p H^n(F';\C).
\]
The desired upper bound follows now from \eqref{eq81}.
\end{proof}

\begin{remark}
\begin{enumerate}
\item Summing \eqref{eq82} over $p$ recovers
$b_n(F)\leq (d-1)^{n+1}= \dim_{\C} M(P_F)$.
\item One can express the upper bound in \eqref{eq82} in terms of binomial coefficients involving only $n$, $d$ and $p$. 
\end{enumerate}
\end{remark}


\section{An example}

\begin{example}
\label{ex1} In this example, we work out the bounds from the previous results on a particular example.

For $n\geq2$, consider the extension of Broughton's example
studied in \cite{MM26}:
\[
 f:\C^{n+1}\rightarrow\C,\quad
 f(x_1,\ldots,x_{n+1})=x_1^2x_2+x_1.
\]
We recall the geometric information needed here.
There are no affine critical points, $B=\{0\}$, and
\[
 F_c\simeq\C^*\times\C^{n-1}\quad(c\ne0),
 \qquad
 F_0\simeq\C^n\sqcup
 \left(\C^*\times\C^{n-1}\right).
\]
The reduced cohomology of the general fiber is
\begin{equation}\label{eq79}
 \widetilde H^k(F_c;\Q)\simeq
 \begin{cases}
  \Q(-1),&k=1,\\
  0,&k\ne1.
 \end{cases}
\end{equation}
For the atypical fiber,
\begin{equation}\label{eq80}
 \widetilde H^k(F_0;\Q)\simeq
 \begin{cases}
  \Q(0),&k=0,\\
  \Q(-1),&k=1,\\
  0,&k\ge2.
 \end{cases}
\end{equation}
These are isomorphisms of mixed Hodge structures.

For the graph completion $\oX$, we obtained in \cite[Example 6.1]{MM26} that 
\[
 \Sigma_0=\Sigma_0^\infty
 =\{x_0=x_1=0,\ t=0\}\simeq\bP^{n-1}.
\]
Its unique top-dimensional stratum is
\[
 S=\Sigma_0^\infty\cap\{x_2\ne0\}\simeq\C^{n-1},
\]
and the transverse calculation in loc.cit. yields
\[
 s=s_0=n-1,\quad
 \mu_S^\pitchfork=\dim_\Q H_S^\pitchfork=1.
\]
We now determine the mixed Hodge structures entering
the bound of our Theorem~\ref{thm3}.

The stalk $\Nu_S$ has dimension one. Also, since $S$ is simply
connected, we get by evaluation
\[
 K_S=H^0(S;\cV_{S,1})\simeq(\Nu_S)_1
\]
as mixed Hodge structures, so $\dim_\Q K_S\leq1$.
As $n-s-1=0$, Theorem~\ref{thm3} yields
the monomorphisms
\[
 \Q(0)\simeq\widetilde H^0(F_0;\Q)
 \hookrightarrow
 \ker\!\left(N_S(1):K_S(1)\to K_S\right)
 \hookrightarrow K_S(1).
\]
The source has dimension one and the final target has
dimension at most one. Both arrows are therefore
isomorphisms, and
\[
 K_S\simeq(\Nu_S)_1=\Nu_S\simeq\Q(-1),
 \quad N_S=0.
\]
In particular, the atypical fiber bound is already an equality
of mixed Hodge structures (not just of dimensions).

By Theorem~\ref{thm2}, since $n-s=1$, we get
\[
 (H_S^\pitchfork)_1^\vee(-2)
 \simeq(\Nu_S)_1\simeq\Q(-1).
\]
Twisting by $(2)$ and dualizing, we obtain
$(H_S^\pitchfork)_1\simeq\Q(-1)$.
Since $\dim_\Q H_S^\pitchfork=1$, we have
$H_S^\pitchfork=(H_S^\pitchfork)_1\simeq\Q(-1)$.
The local Hodge bound is therefore an equality:
\[
 h^{p,q}\left(\widetilde H^0(F_0;\Q)\right)
 =
 h_1^{1-p,1-q}(H_S^\pitchfork)
 =
 \begin{cases}
  1,&(p,q)=(0,0),\\
  0,&\text{otherwise}.
 \end{cases}
\]

We can also use this example to distinguish nonzero vanishing
cohomology from nontrivial monodromy on the general fiber.
Indeed, the algebraic isomorphism
\[
 F_1\times\C^*\xrightarrow{\sim}f^{-1}(\C^*),
 \qquad
 ((x_1,x_2,x_3,\ldots,x_{n+1}),c)
 \longmapsto
 (cx_1,c^{-1}x_2,x_3,\ldots,x_{n+1})
\]
identifies $f$ with projection to $c$. Thus its cohomology
local systems on $\C^*$ are constant, and
$T_0=\operatorname{id}$ on $H^1(F_c;\Q)$.

On the other hand, since $B=\{0\}$, the isomorphism
$V^1(0)\simeq\widetilde H^1(F_c;\Q)$ from \cite{MM26, Br03} implies
$\dim_\Q V^1(0)=1$. The injection
\[
 V_H^1(0)\hookrightarrow H^0(S;\cV_S)
 =K_S\simeq\Q(-1)
\]
is therefore an isomorphism of mixed Hodge structures.
In particular,
\[
 V_H^1(0)\simeq\Q(-1),
 \quad {\rm but} \ \  
 \operatorname{rk}(T_0-\operatorname{id})
 =0<\dim_\Q V^1(0)=1.
\]
So the Hodge bound for the atypical fiber is sharp, while the
monodromy rank bound \eqref{eq43}
is strict.\hfill$\qed$
\end{example}



\begin{thebibliography}{10}


\bibitem{Br03}
T.~Br\'elivet,
\emph{Topologie des fonctions r\'eguli\`eres et cycles
\'evanescents},
Rev. Mat. Complut. {16} (2003), no.~1, 131--149.

\bibitem{Br88}
S.~A. Broughton,
\emph{Milnor numbers and the topology of polynomial hypersurfaces},
Invent. Math. {92} (1988), no.~2, 217--241.

\bibitem{D71}
P.~Deligne,
\emph{Th\'eorie de Hodge. II},
Publ. Math. Inst. Hautes \'Etudes Sci. {40} (1971), 5--57.

\bibitem{D74}
P.~Deligne,
\emph{Th\'eorie de Hodge. III},
Publ. Math. Inst. Hautes \'Etudes Sci. {44} (1974), 5--77.

\bibitem{D04}
A.~Dimca,
\emph{Sheaves in Topology},
Universitext, Springer-Verlag, Berlin, 2004.

\bibitem{DN01}
A.~Dimca and A.~N\'emethi,
\emph{On the monodromy of complex polynomials},
Duke Math. J. {108} (2001), no.~2, 199--209.

\bibitem{DS93}
A.~Dimca and M.~Saito,
\emph{On the cohomology of a general fiber of a polynomial map},
Compositio Math. {85} (1993), no.~3, 299--309.

\bibitem{DS03}
A.~Dimca and M.~Saito,
\emph{Monodromy at infinity and the weights of cohomology},
Compositio Math. {138} (2003), no.~1, 55--71.


\bibitem{GLN96}
R.~Garc\'ia L\'opez and A.~N\'emethi,
\emph{On the monodromy at infinity of a polynomial map},
Compositio Math. {100} (1996), no.~2, 205--231;
appendix by R.~Garc\'ia L\'opez and J.~H.~M. Steenbrink.

\bibitem{GLN99a}
R.~Garc\'ia L\'opez and A.~N\'emethi,
\emph{On the monodromy at infinity of a polynomial map, II},
Compositio Math. {115} (1999), no.~1, 1--20.

\bibitem{GLN99b}
R.~Garc\'ia L\'opez and A.~N\'emethi,
\emph{Hodge numbers attached to a polynomial map},
Ann. Inst. Fourier (Grenoble) {49} (1999), no.~5,
1547--1579.




\bibitem{MT13}
Y.~Matsui and K.~Takeuchi,
\emph{Monodromy at infinity of polynomial maps and Newton
polyhedra}, with an appendix by C.~Sabbah,
Int. Math. Res. Not. IMRN (2013), no.~8, 1691--1746.

\bibitem{M19} L. Maxim, \emph{Intersection homology \& perverse sheaves--with applications to singularities}, Graduate Texts in Mathematics, 281. Springer, Cham, 2019.

\bibitem{MM26}
L.~Maxim and J.~Messina,
\emph{On the topology of fibers of complex polynomial maps},
preprint arXiv:2608.04320.

\bibitem{MPT22}
L.~Maxim, L.~P\u{a}unescu, and M.~Tib\u{a}r,
\emph{The vanishing cohomology of non-isolated hypersurface
singularities},
J. Lond. Math. Soc. (2) {106} (2022), no.~1, 112--153.

\bibitem{MSS20}
L.~Maxim, M.~Saito, and J.~Sch\"urmann,
\emph{Spectral Hirzebruch--Milnor classes of singular hypersurfaces},
Math. Ann. {377} (2020), no.~1--2, 281--315.

\bibitem{MS22}
L.~Maxim and J.~Sch\"urmann,
\emph{Constructible sheaf complexes in complex geometry and
applications},
in \emph{Handbook of Geometry and Topology of Singularities III},
Springer, Cham, 2022, 679--791.

\bibitem{NS99}
A.~N\'emethi and C.~Sabbah,
\emph{Semicontinuity of the spectrum at infinity},
Abh. Math. Sem. Univ. Hamburg {69} (1999), 25--35.

\bibitem{NN00}
W.~D. Neumann and P.~Norbury,
\emph{Vanishing cycles and monodromy of complex polynomials},
Duke Math. J. {101} (2000), no.~3, 487--497.

\bibitem{PS08} C..~Peters and J.~Steenbrink,
\emph{Mixed Hodge Structures}, Ergeb. Math. Grenzgeb. (3), vol.~52, Springer, Berlin, 2008.

\bibitem{Sab97}
C.~Sabbah,
\emph{Monodromy at infinity and Fourier transform},
Publ. Res. Inst. Math. Sci. {33} (1997), no.~4, 643--685.

\bibitem{Sa88}
M.~Saito,
\emph{Modules de Hodge polarisables},
Publ. Res. Inst. Math. Sci. {24} (1988), no.~6, 849--995.

\bibitem{Sa89}
M.~Saito,
\emph{Duality for vanishing cycle functors},
Publ. Res. Inst. Math. Sci. {25} (1989), no.~6, 889--921.

\bibitem{Sa90}
M.~Saito,
\emph{Mixed Hodge modules},
Publ. Res. Inst. Math. Sci. {26} (1990), no.~2, 221--333.

\bibitem{SS85}
J.~Scherk and J.~H.~M. Steenbrink,
\emph{On the mixed Hodge structure on the cohomology of the Milnor
fibre},
Math. Ann. {271} (1985), no.~4, 641--665.

\bibitem{Sch73}
W.~Schmid,
\emph{Variation of Hodge structure: the singularities of the period
mapping},
Invent. Math. {22} (1973), 211--319.

\bibitem{Schu03}
J.~Sch\"urmann,
\emph{Topology of Singular Spaces and Constructible Sheaves},
Monografie Matematyczne (New Series), vol.~63,
Birkh\"auser Verlag, Basel, 2003.

\bibitem{ST95}
D.~Siersma and M.~Tib\u{a}r,
\emph{Singularities at infinity and their vanishing cycles},
Duke Math. J. {80} (1995), no.~3, 771--783.

\bibitem{ST00}
D.~Siersma and M.~Tib\u{a}r,
\emph{Singularities at infinity and their vanishing cycles. II.
Monodromy},
Publ. Res. Inst. Math. Sci. {36} (2000), no.~5, 659--679.



\bibitem{ST07}
D.~Siersma and M.~Tib\u{a}r,
\emph{Singularity exchange at the frontier of the space}, Real and complex singularities, 327--342, Trends Math., Birkh\"auser, Basel, 2007.

\bibitem{ST11} D.~Siersma and M.~Tib\u{a}r,  \emph{Betti bounds of polynomials}, 
Mosc. Math. J. 11 (2011), no. 3, 599--615, 632.

\bibitem{Ste76}
J.~H.~M. Steenbrink,
\emph{Limits of Hodge structures},
Invent. Math. {31} (1976), no.~3, 229--257.

\bibitem{Ste77} J.~H.~M. Steenbrink,
\emph{Intersection form for quasi-homogeneous singularities},
Compos. Math. 34 (1977), 211--223.

\bibitem{SZ85}
J.~H.~M. Steenbrink and S.~Zucker,
\emph{Variation of mixed Hodge structure. I},
Invent. Math. {80} (1985), no.~3, 489--542.

\bibitem{TT16}
K.~Takeuchi and M.~Tib\u{a}r,
\emph{Monodromies at infinity of non-tame polynomials},
Bull. Soc. Math. France {144} (2016), no.~3, 477--506.

\bibitem{T07}
M.~Tib\u{a}r,
\emph{Polynomials and Vanishing Cycles},
Cambridge Tracts in Mathematics, vol.~170,
Cambridge University Press, Cambridge, 2007.


\end{thebibliography}
\end{document}